\documentclass{amsart}

\usepackage{amsmath,amsthm,amssymb}

\renewcommand{\frak}{\mathfrak}
\renewcommand{\geq}{\geqslant}
\renewcommand{\leq}{\leqslant}
\newcommand{\bbold}{\mathbb}
\newcommand{\rom}{\textup}
\renewcommand\emptyset{\varnothing}
\newcommand{\abs}[1]{\lvert#1\rvert}

\def \R {{\bbold R}}
\def \Q {{\bbold Q}}
\def \Z {{\bbold Z}}
\def \N {{\bbold N}}
\def \C {{\bbold C}}
\def \F {{\bbold F}}
\def \m {{\frak m}}
 
\def \bar {\overline}
\def \<{\langle}
\def \>{\rangle}

\def \iso {\cong}
\def \hat{\widehat}
\def \rank{\operatorname{rank}}
 
\def \OX {{\cal O}\langle X\rangle}
\def \OhatX {\hat{\cal O}\langle X\rangle}
\def \OhatXprime {\hat{\cal O}\langle X'\rangle}

\def \OXprime {{\cal O}\langle X'\rangle}
\def \ZpX {\Z_p\langle X\rangle}

\def \Frac {\operatorname{Frac}}

\def \Sol {\operatorname{Sol}}
\def \cal {\mathcal}
\def \tr {\operatorname{tr}}

\newtheorem{theorem}{Theorem}[section]
\newtheorem*{theoremIntro}{Theorem}
\newtheorem*{theoremA}{Theorem A}
\newtheorem*{theoremB}{Theorem B}
\newtheorem{lemma}[theorem]{Lemma}
\newtheorem{prop}[theorem]{Proposition}
\newtheorem{cor}[theorem]{Corollary}

\theoremstyle{definition}

\theoremstyle{remark}
\newtheorem{remarkNumbered}[theorem]{Remark}

\newtheorem*{example}{Example}

\newtheorem*{remarks}{Remarks}
\newtheorem*{remark}{Remark}

\numberwithin{equation}{section}

\begin{document}

\title{Ideal Membership in Polynomial Rings over the Integers}

\author{Matthias Aschenbrenner}
\address{Mathematical Sciences Research Institute \\
17 Gauss Way\\
Berkeley, CA 94720, U.S.A. }

\address{Department of Mathematics\\
University of California at Berkeley\\
Evans Hall \\
Berkeley, CA 94720, U.S.A. }
\email{maschenb@math.berkeley.edu}

\date{June 2003}

\begin{abstract} 
We present a new approach to the ideal membership problem for polynomial
rings over the integers: given polynomials $f_0,f_1,\dots,f_n\in\Z[X]$,
where $X=(X_1,\dots,X_N)$ is an $N$-tuple of indeterminates, 
are there $g_1,\dots,g_n\in\Z[X]$ such that
$f_0=g_1f_1+\cdots+g_nf_n$? We show that the degree of the polynomials
$g_1,\dots,g_n$ can be bounded by $(2d)^{2^{O(N^2)}}(h+1)$ where
$d$ is the maximum total degree and $h$ the maximum height of the
coefficients of $f_0,\dots,f_n$.
Some related questions, primarily 
concerning linear equations in $R[X]$, where $R$
is the ring of integers of a number field,
are also treated.
\end{abstract}


\maketitle

\section*{Introduction}

The following well-known theorem, due to Grete Hermann \cite{Hermann}, 1926,
gives an upper bound on the complexity of the ideal membership problem
for polynomial rings over fields:

\begin{theoremIntro}
Consider polynomials $f_0,\dots,f_n\in F[X]=F[X_1,\dots,X_N]$ 
of \rom{(}total\rom{)} degree $\leq d$ over
a field $F$. If $f_0\in (f_1,\dots,f_n)$, then
$$f_0=g_1f_1+\cdots+g_nf_n$$
for certain $g_1,\dots,g_n\in F[X]$ whose degrees are bounded by 
$\beta$, where $\beta=\beta(N,d)$ depends only on $N$ and $d$
\rom{(}and not on the field $F$ or the particular polynomials $f_0,\dots,f_n$\rom{)}.
\end{theoremIntro}

This theorem was a first step in Hermann's  project, initiated 
by work of Hentzelt and Noether \cite{Hentzelt-Noether}, to construct bounds
for some of the central operations of commutative algebra in polynomials rings
over fields. A simplified and corrected proof was published by Seidenberg 
\cite{Seidenberg1} in the 1970s, with an explicit but incorrect 
bound $\beta(N,d)$ (see \cite{Renschuch}). We will reproduce 
a proof, using Hermann's classical method, in Section~\ref{HermannSection} 
below and show that one may take
$$\beta(N,d)=(2d)^{2^{N}}.$$
Note that the computable character of this bound
reduces the question whether
$f_0\in (f_1,\dots,f_n)$ for given $f_j\in F[X]$ to solving an (enormous)
system of linear equations over $F$. Hence, in this way one obtains a
(naive) algorithm to solve the ideal membership problem for $F[X]$
(provided $F$ is given in some explicitly computable manner).
Later, Buchberger in his Ph.D.~thesis (1965)
introduced the important concept of a {\it Gr\"obner basis,}\/
and gave an algorithm for deciding ideal membership for $F[X]$, which is
widely used today (see, e.g., \cite{Becker-Weispfenning}). 

The doubly exponential
nature of $\beta$ above is essentially unavoidable, 
as a family of examples
due to Mayr and Meyer \cite{Mayr-Meyer} show. 
In fact, they prove that ideal membership for $\Q[X]$ is
exponential-space hard:
the amount of space needed by {\it any}\/
algorithm to decide ideal membership for $\Q[X]$ (or $\Z[X]$) grows
exponentially in the size of the input.
If we restrict to $f_0,\dots,f_n$ of a special
form, often dramatic improvements are possible: 
for example, if $f_0=1$ (the situation of Hilbert's Nullstellensatz), 
then in the theorem we may replace the doubly exponential 
$(2d)^{2^{N}}$ by the single exponential bound $d^N$, if $d>2$ (due to
Koll\'ar \cite{Kollar}), and by $2^{N+1}$ if $d=2$ (due to Sombra
\cite{Sombra}). 
A number of results
show the existence of single-exponential bounds in the (general) ideal 
membership problem for $F[X]$, 
under suitable geometric assumptions on the ideal $I=(f_1,\dots,f_n)$: 
for example if $I$ is zero-dimensional or
a complete intersection \cite{Berenstein-Yger}, or unmixed \cite{DFGS}.

In this paper, we study the ideal membership problem 
over coefficient rings of an arithmetic nature, like the ring of integers $\Z$
(instead of over a field $F$).
An easy example (see Section~\ref{Ideal-Section}) 
shows that contrary to what happens over fields, if a bound $d$ on the degree of
$f_0,f_1,\dots,f_n\in\Z[X]$ is given and $f_0\in (f_1,\dots,f_n)\Z[X]$,
then {\it there is no uniform bound on the degrees of $g_j$'s such that
$f_0=g_1f_1+\cdots+g_nf_n$, which depends only on $N$ and $d$.}\/ 
So any bound on the degree of the $g_1,\dots,g_n$
as a function of $f_0,f_1,\dots,f_n$ will necessarily also have to depend on the coefficients
of the polynomials $f_j$.

A decision procedure for
the ideal membership problem for polynomial rings over $\Z$ 
has been known at least since the early 1970s, see, e.g.,
\cite{Ayoub}, \cite{Baur}, \cite{Evans},  \cite{KRK},
\cite{Richman}, \cite{Seidenberg3}, \cite{Simmons}.
However, these results did not yield the existence of a {\it primitive}\/ 
recursive algorithm, for any fixed $N\geq 3$, let alone the existence of 
bounds similar to the ones in Hermann's theorem for polynomial rings over
fields. Indeed, it was suspected by some that this was one of the rare
cases where a natural decision problem allows an algorithmic solution,
but not a primitive recursive one.  (See
\cite{kron} for a survey of the history and the various proposals for
computing in $\Z[X]$.)

Finding a decision procedure for ideal membership in $\Z[X]$ was 
central to Kronecker's ideology of constructive mathematics \cite{Edwards}. 
In fact, one may argue that he was primarily interested in 
what we would call today a
{\em primitive recursive}\/ algorithm. 
Thus, the task of finding a primitive recursive decision method for
ideal membership in $\Z[X]$
has aptly been called ``Kronecker's problem'' in \cite{Gallo-Mishra}.
In this paper, Gallo and Mishra 
adapted Buchberger's algorithm for the construction of Gr\"obner bases and deduced a primitive
recursive procedure to decide the ideal membership problem for $\Z[X]$, when the number
of variables $N$ is fixed. Analyzing their algorithm, they obtained the following bounds:

\begin{theoremIntro}
Let $f_0,\dots,f_n\in\Z[X]$. If $f_0 \in (f_1,\dots,f_n)$, then
$$f_0=g_1f_1+\cdots+g_nf_n$$
for certain polynomials $g_1,\dots,g_n\in\Z[X]$ whose size $|g_j|$ is bounded by 
$$W_{4N+8}\bigl(|f_0|+\cdots+|f_n|+N\bigr).$$
\end{theoremIntro}

Here the size $|f|$ of a polynomial $f\in\Z[X]$ is a crude measure of its 
complexity, and equals the
maximum of the absolute values of the coefficients and the 
degrees of $f$ with respect to each indeterminate,
see \cite{Gallo-Mishra}, p.~346. The function $W_k$ is the $k$th function in the so-called Wainer
hierarchy of primitive recursive functions, see \cite{Wainer}. Even for small $k$,  
these functions are already very rapidly growing:  We have $W_0(n)=n+1$, $W_1(n)=2n+1$, but
$W_2$ grows asymptotically like the exponential $n\mapsto 2^n$, and $W_3$ as the 
$n$-times iterated exponential function, and so on. These bounds are only primitive recursive for each
fixed $N$, the growth rate of this bound as a function of $N$ being similar to
the notorious Ackermann function. 

Gallo and Mishra's analysis of the complexity of
their algorithm ultimately rests on an
effective version of
Hilbert's Basis Theorem for increasing chains of monomial ideals in $\Z[X]$.
This approach is doomed to fail in providing bounds which 
are also primitive recursive 
for varying $N$: In general, 
even the length of an increasing 
chain of ideals in $\Z[X]$ 
with the $k$-th ideal in the chain generated by monomials of degree $\leq kd$,
can have a  growth behavior similar to Ackermann's function,
as a function of $N$ and $d$.
(See \cite{Moreno-Socias}.) 

In the present paper, we will give a proof of the following theorem. Given
a polynomial $f\in\Z[X]$ 
we let $h(f)$ be the height of $f$, that is, the maximum of $\log\,\abs{a}$
where $a$ ranges over the non-zero coefficients of $f$, with $h(0):=0$.

\begin{theoremA}
If $f_0,f_1,\dots,f_n\in\Z[X]$ are polynomials with $f_0\in (f_1,\dots,f_n)$, 
whose degrees are at most $d$ and whose heights are at most $h$,  then
$$f_0=g_1f_1+\cdots+g_nf_n$$
for certain polynomials $g_1,\dots,g_n\in\Z[X]$ of degrees at most
$$\gamma(N,d,h)=(2d)^{2^{O(N^2)}}(h+1).$$
\end{theoremA}

In principle, the (universal) constant hidden in the $O$-notation can 
be made explicit,
see Section~\ref{Ideal-Section} below.
The bound $\gamma$ on the degrees of the $g_j$'s implies the existence
of a similar (doubly exponential) bound on the heights of the $g_j$.
As a consequence, we obtain a naive elementary recursive
decision procedure for ideal membership in $\Z[X]$.
In this paper we prove in fact a generalization of
Theorem~A with $\Z$ replaced by the ring of integers of a
number field $F$, using an appropriate
notion of height for elements of $F$.

The starting point for our proof of Theorem~A 
is the simple observation that one can {\em localize}\/
the question whether $f_0\in (f_1,\dots,f_n)$ and reduce it to
finitely many subproblems in the following way:
Using (and refining) the
classical method of Hermann (for $F=\Q$), one can test whether 
$f_0\in (f_1,\dots,f_n)\Q[X]$, and assuming this is so, we obtain, by clearing denominators, a
representation
\begin{equation}\label{1}\tag{1}
\delta f_0 = g_1f_1 + \cdots + g_nf_n \qquad\text{with $\delta\in\Z$, 
$\delta\neq 0$, $g_1,\dots,g_n\in\Z[X]$.}
\end{equation}
Let $p_1,\dots,p_K$ be the different prime factors of $\delta$. Then another necessary
condition for $f_0\in (f_1,\dots,f_n)$, besides $f_0\in (f_1,\dots,f_n)\Q[X]$, is that
$f_0\in (f_1,\dots,f_n)\Z_{(p_k)}[X]$ for $k=1,\dots,K$. Together with \eqref{1}, these 
{\em necessary}\/ conditions are also {\em sufficient}\/ for $f_0\in (f_1,\dots,f_n)$: If
$f_0\in (f_1,\dots,f_n)\Z_{(p_k)}[X]$, then
\begin{equation}\label{2}\tag{2}
\delta_k f_0 = g_{1k} f_1 + \cdots + g_{nk} f_n \qquad \text{for some $\delta_k\in\Z\setminus p_k\Z$ and
$g_{jk}\in\Z[X]$.}
\end{equation}
Since $\delta,\delta_1,\dots,\delta_K$ have no common prime factor we find, by the Euclidean Algorithm,
a linear combination of them that equals $1$:
\begin{equation}\label{3}\tag{3}
a\delta + a_1\delta_1 + \cdots + a_K\delta_K=1 \qquad (a,a_1,\dots,a_K\in\Z).
\end{equation}
Combining \eqref{1}, \eqref{2} and \eqref{3} we get
$$f_0=\bigl(a\delta+a_1\delta_1+\cdots+a_K\delta_K\bigr)f_0=
\sum_{j=1}^n \bigl(ag_j+a_1g_{j1}+\cdots+a_Kg_{jK}\bigr)f_j,$$
which exhibits $f_0$ as an element of $(f_1,\dots,f_n)$.

Now note that given a prime $p$, we have
$f_0\in (f_1,\dots,f_n)\Z_{(p)}[X]$ if and only if
the homogeneous linear equation
\begin{equation}\label{4}\tag{4}
f_1y_1+\cdots+f_ny_n-f_0y_{n+1}=0
\end{equation}
in the unknowns $y_1,\dots,y_{n+1}$ has a solution $(y_1,\dots,y_{n+1})\in
\bigl(\Z_{(p)}[X]\bigr)^{n+1}$ with $y_{n+1}=1$. This reduces the
problem of deciding whether $f_0\in (f_1,\dots,f_n)\Z_{(p)}[X]$ to the
two following subproblems:
\begin{itemize}
\item[(a)] constructing a collection of generators
$z^{(1)},\dots,z^{(L)}\in\bigl(\Z_{(p)}[X]\bigr)^{n+1}$ for the module
of solutions (in $\Z_{(p)}[X]$) to the equation \eqref{4}, and
\item[(b)] deciding whether the ideal in $\Z_{(p)}[X]$ generated by the
last components of the vectors $z^{(1)},\dots,z^{(L)}$ contains $1$.
\end{itemize}
Problem~(b) can be easily treated by applying 
the effective Nullstellensatz for
$\Q[X]$ and $\F_p[X]$ (or Hermann's Theorem). 
By a faithful flatness argument, it is possible
to further reduce problem (a) to the construction of a set 
of generators for the $\Q[X]$-module of solutions to \eqref{4} in
$\Q[X]$, and a set of generators $S\subseteq\bigl(\Z_{(p)}[X]\bigr)^{n+1}$ for 
the $\ZpX$-module of solutions to \eqref{4} in $\ZpX$.
Here, $\ZpX$ denotes the ring of restricted power
series with $p$-adic integer coefficients (see ~\cite{BGR} or 
Section~\ref{PowerSeriesSection}). 
The great advantage of the power series rings $\ZpX$ over polynomial rings over $\Z$ (or over $\Z_{(p)}$) is
that they satisfy a Weierstra\ss{} Division and Preparation Theorem.
Hermann's method for deciding ideal membership, that is, deciding
solvability of a single inhomogeneous linear equation, has a variant
which allows for the
construction of a finite set of generators for the $\Q[X]$-module of
solutions to the linear homogeneous equation \eqref{4} in $\Q[X]$.
The key step in our argument is to adapt this method to 
explicitly construct the set $S$ from above, that is, to show
the {\em effective flatness}\/ of $\ZpX$ as $\Z_{(p)}[X]$-module.
All computations take place in $\Z_{(p)}[X]$, and
bounds for the heights of the polynomials
occurring in each step can be found. This enables us to
calculate the bound $\gamma$.

Theorem~A naturally generalizes 
to {\em systems of linear equations}\/ over polynomial rings, and as
the sketch above already indicates, one also obtains 
information on homogeneous systems of linear equations. 
For example, the methods developed here lead to the following theorem
on degree bounds for generators of syzygies:

\begin{theoremB}
The $\Z[X]$-module of solutions $(y_1,\dots,y_n)\in\bigl(\Z[X]\bigr)^n$ 
of the equation
$$f_1y_1+\cdots+f_ny_n=0,$$ where $f_1,\dots,f_n\in\Z[X]$ are of degree
$\leq d$, is generated by solutions 
$$y^{(1)},\dots,y^{(K)}\in\bigl(\Z[X]\bigr)^n$$
whose entries are of degree $\leq (2d)^{2^{O(N^2)}}$.
\end{theoremB}

Note that 
this bound does not depend on the coefficients of the $f_j$'s.
(The number $K$ of required generators also depends only on
$N,n$ and $d$.) Theorem~B holds in a rather more general context, 
for any almost hereditary ring in place of $\Z$. 
See Section~\ref{Effective-Flatness-Section} below for
the definition of almost hereditary rings, and
\cite{coherent} on uniform degree bounds on syzygies for an even 
larger class of rings.
The size of the coefficients of the entries of $y^{(k)}$ can be similarly
estimated, by a bound that also depends on the heights
of $f_1,\dots,f_n$. 

\subsection*{Organization of the paper}
We begin (in Section~\ref{Abs-Value-Section}) by recalling basic
definitions about absolute values on number fields, defining a height
function on the algebraic closure of $\Q$, and establishing some auxiliary
facts about it used later. In Section~\ref{PowerSeriesSection} we state 
some fundamental facts about the ring of restricted power 
series over a complete discrete valuation ring. Section~\ref{HermannSection}
contains an exposition of Hermann's method for solving  
systems of linear equations in polynomial rings over fields.
This is the basis for Section~\ref{Effective-Flatness-Section}, 
where we give a proof of Theorem~B modeled on this method. 
We also indicate two applications
concerning bounds for some operations on finitely generated modules and 
a criterion for primeness of ideals in $\Z[X]$ in the style
of \cite{Schmidt-Goettsch}. In Section~\ref{Height-Section} we complement
Theorem~B for rings of integers in number fields by establishing bounds on
the height of generators for syzygy modules. In 
Section~\ref{Ideal-Section}
we use these results to prove Theorem~A.

\subsection*{Acknowledgments}
This paper is based on a part of the author's 
Ph.D.~thesis \cite{Aschenbrenner-thesis} written
under the direction of Lou van den Dries, whom he would like to 
thank for his guidance and advice.
He would also like to thank Hendrik Lenstra and
Bjorn Poonen for useful suggestions
concerning the proof of Lemma~\ref{C-Lemma}.

\subsection*{Notations and conventions}
Throughout this paper
$\N=\{0,1,2,\dots\}$
denotes the set of natural numbers. 

Let $R$ be a ring (here and below: always commutative with a unit element). 
The localization $S^{-1}R$, where $S$ denotes the set of non-zero-divisors
of $R$, is called the ring of fractions of $R$, denoted by $\Frac(R)$.
If $A$ is an $m\times n$-matrix with entries in $R$, the set of solutions
in $R^n$ to the homogeneous system of linear equations $Ay=0$ 
is an $R$-submodule of $R^n$, which we denote by $\Sol_R(A)$.
It is sometimes called the (first) {\bf module of syzygies of $A$.} If
$R$ is coherent (e.g., if $R$ is Noetherian), 
then $\Sol_R(A)$ is finitely generated.
For submodules $M$, $M'$  of an $R$-module we write
$$(M': M):=\bigl\{a\in R:am\in M' \text{ for all $m\in M$}\bigr\},$$
an ideal of $R$ (containing the annihilator of $M$).

By $X=(X_1,\dots,X_N)$ we always denote a tuple of $N$ distinct 
indeterminates, where $N\in\N$. The (total)
degree of a polynomial $0\neq f\in R[X]=R[X_1,\dots,X_N]$ is denoted
by $\deg(f)$, and
the degree of $f$ in $X_i$ (where $i\in\{1,\dots,N\}$) by $\deg_{X_i}(f)$.
By convention $\deg(0) := -\infty$ and $\deg_{X_i}(0) := -\infty$,
where $-\infty<\N$.
We extend this notation to finite tuples
$f=(f_1,\dots,f_n)$ of polynomials in $R[X]$ by setting
$\deg(f):=\max_j \deg(f_j)$ (the degree of $f$).
Similarly we define $\deg_{X_i}(f)$.

The notions of {\it computable field}\/ and {\it computable ring}\/ 
are used in an
informal way. 
We will say that a computable ring $R$ is {\bf syzygy-solvable} if
there is an algorithm which, given $a_1,\dots,a_n\in R$, constructs
a finite set of generators for the solutions to the homogeneous
linear equation $a_1y_1+\cdots+a_ny_n=0$. 
(This is called ``finitely related'' in \cite{Richman}.) 
For example, the prime fields
$\Q$ and $\F_p$  are clearly syzygy-solvable, as is
$\Z$, or more generally the ring of integers of any number field
(see \cite{Cohen}).

\section{Absolute Values and Height Functions}\label{Abs-Value-Section}

We assume that the reader is familiar with the basic theory of absolute
values on number fields as expounded in, say, \cite{Lang}, Chapter~II,
and the (absolute, logarithmic) height function  on the algebraic closure
of $\Q$ as used in diophantine geometry (see \cite{Lang-Diophantine}, 
Chapter~3).
We recall some definitions and a few basic facts used later. 

\subsection*{Absolute values}
We let $\abs{\,\cdot\,}$ denote the usual (Euclidean) absolute value 
$\abs{\,\cdot\,}$ on $\Q$, and
for a prime number $p$ we let $\abs{\,\cdot\,}_p$ denote the $p$-adic
absolute value on $\Q$: $\abs{a}_v=p^{-v_p(a)}$ for 
$a\in\Q^\times$, where $v_p\colon\Q^\times\to\Z$ denotes the $p$-adic
valuation on $\Q$.
Let $F$ be an algebraic number field of degree 
$d=[F:\Q]$, and let $R={\cal O}_F$ be the ring of integers of $F$. 
If $v$ is a finite place of $F$ which lies over the
prime number $p$, we write $v|p$. If $v$ is an infinite place of $F$ we
write $v|\infty$. 
To every place $v$ of $F$ we associate an absolute value 
$\abs{\,\cdot\,}_v$ on $F$, normalized so that
\begin{enumerate}
\item if $v|p$ for a prime $p$, then $\abs{\,\cdot\,}_v$ extends the $p$-adic
absolute value  $\abs{\,\cdot\,}_p$ on $\Q$,
\item if $v|\infty$, then $\abs{\,\cdot\,}_v$ extends the usual absolute
value $\abs{\,\cdot\,}$ on $\Q$.
\end{enumerate}
For every place $v$ of $F$ we let $F_v$ denote the
completion of $F$ with respect to the topology induced by $\abs{\,\cdot\,}_v$ 
and $\Q_v\subseteq F_v$ the completion
of $\Q$ with respect to the topology induced by the 
restriction of $\abs{\,\cdot\,}_v$ to $\Q$. We put 
$d_v = [F_v:\Q_v]$. If $v|p$ is finite, then 
$\Q_p$ is the field of $p$-adic numbers. 
If $v|\infty$, then $\Q_v=\R$, and
either $F_v=\R$ and $d_v=1$ (in which case $v$ is called real), or 
$F_v=\C$ and $d_v=2$ ($v$ is complex).
Given $w=\infty$ or $w=p$ for a prime $p$ we have 
$$d=\sum_{v|w} d_v,$$
where the sum ranges over all places $v$ of $F$ with $v|w$.
We let $M_F$ denote the set of all places of $F$, 
$M_F^\infty:=\{v\in M_F: v|\infty\}$ the set of
infinite places, and
$M_F^0 := M_F\setminus M_F^\infty$ the set of finite places of $F$.
The number field $F$ satisfies
the following {\bf product formula} (with multiplicities $d_v$):
\begin{equation}\label{Product}
\prod_{v\in M_F} ||a||_v = 1 \qquad (a\in F^\times),
\end{equation}
where $||a||_v:=\abs{a}_v^{d_v}$
for $v\in M_F$.

The assignment $v\mapsto
{\frak p}_v:=\bigl\{r\in R:\abs{r}_v<1\bigr\}$ establishes a one-to-one
correspondence between $M_F^0$ and the set of non-zero prime ideals of $R$.
If $v|p$ is a 
finite place of $F$ and ${\frak p}={\frak p}_v$,
then the absolute value $\abs{\,\cdot\,}_v$ on $F$ associated to $v$ and the
${\frak p}$-adic valuation on $F$ are connected as follows:
$$\abs{a}_v=p^{-v_{\frak p}(a)/e_v}\qquad\text{for all $a\in F^\times$.}$$
Here $e_v$ denotes the
ramification index of $v$, that is, the unique integer such that 
$p=\pi^{e_v}u$ for some unit $u$ of $R_{\frak p}$ and some
$\pi\in R_{\frak p}$ with $v_{\frak p}(\pi)=1$. We have $e_v|d_v$, 
in fact, $\#(R/\frak p)=p^{d_v/e_v}$.

%

\subsection*{$M_F$-divisors}
An {\bf $M_F$-divisor} is a function ${\frak c}\colon M_F\to\R$ such that
\begin{enumerate}
\item ${\frak c}(v)>0$ for all $v\in M_F$;
\item ${\frak c}(v)=1$ for all but finitely many $v\in M_F$;
\item if $v\in M_F^0$, then there exists an element $a\in F$ with
${\frak c}(v)=\abs{a}_v$.
\end{enumerate}
We shall sometimes write $\abs{\frak c}_v$ instead of ${\frak c}(v)$, and
we put $||\frak c||_v:=\abs{\frak c}_v^{d_v}$. We define the {\bf size}
of an $M_F$-divisor $\frak c$ to be $$||\frak c||_F := \prod_v ||\frak c||_v.$$
The product $\frak c\cdot\frak d$ of two $M_F$-divisors $\frak c$ and $\frak d$
is an $M_F$-divisor, and $||\frak c\cdot\frak d||_F=||\frak c||_F\cdot 
||\frak d||_F$. Given an $M_F$-divisor $\frak c$ we let 
$$L(\frak c)=\bigl\{ a\in F: \text{$\abs{a}_v \leq {\frak c}(v)$ for all
$v\in M_F$}\bigr\},$$
a finite set.
Each non-zero fractional ideal $I$ of $F$ (i.e., a finitely generated 
$R$-submodule of $F$)
determines a unique
$M_F$-divisor ${\frak c}_I$ such that $L(\frak c)=I$ and $\frak c(v)=1$ for
all $v\in M_F^\infty$. We have
$${\frak c}_I(v_{\frak p})=p^{-v_{\frak p}(I)/e_{v_{\frak p}}}$$
for all primes $\frak p\neq 0$ of $R$.
Here $I=\prod_{\frak p} {\frak p}^{v_{\frak p}(I)}$ is the unique
representation of $I$ as a product of non-zero prime ideals of $R$, with
$$v_{\frak p}(I) = \min\bigl\{ v_{\frak p}(a) : a\in I\bigr\} \in \Z$$
and $v_{\frak p}(I)=0$ for almost all $\frak p$. We have $||{\frak c}_I||_F=
1/N(I)$, where $$N(I)=\prod_{\frak p} \#(R/\frak p)^{v_{\frak p}(I)}$$
denotes the norm of $I$. If $I\subseteq R$, 
then $N(I)=\#(R/I)$.

\subsection*{Height functions}
Given a place $v\in M_F$ and a non-empty finite set $S\subseteq F$ we put
$\abs{S}_v := \max\bigl\{\abs{a}_v:a\in S\bigr\}$ and $||S||_v:=\abs{S}_v^{d_v}$.
We define the
(logarithmic) {\em local height}\/ $h_v(S)$ of $S$ at $v$ by
$$h_v(S):=\log^+ ||S||_v.$$
Here $\log^+ r:=\max\bigl\{0,\log r\bigr\}$ for $r\in\R^{>0}$ and
$\log^+ 0:=0$. We declare $h_v(\emptyset):=0$. For a
polynomial $f\in F[X]$
we put $||f||_v:=||S||_v$, where $S$ is the
set of coefficients of $f$. The {\em local height}\/ of $f$
at $v\in M_F$ is defined by $$h_v(f):=\log^+ ||f||_v.$$
More generally, for $f_1,\dots,f_n\in F[X]$ we put
$$h_v(f_1,\dots,f_n):=\log^+ ||S||_v,$$ where $S$ is the set of coefficients
of $f_1,\dots,f_n$. Note $h_v(f_1,\dots,f_n)\geq 0$.
Here are some other basic
properties of $h_v$, immediate from the definition:

\begin{lemma}\label{hv-lemma}
Let $v\in M_F$ and $a,a_1,\dots,a_n\in F$. Then:
\begin{enumerate}
\item $h_v(a)=h_v(-a)$;
\item $h_v(a^k)=k\cdot h_v(a)$ for $k\in\N$;
\item $h_v(a_1+\cdots+a_n) \leq h_v(a_1,\dots,a_n) + \log n$ if $v\in M_F^\infty$;
\item $h_v(a_1+\cdots+a_n) \leq h_v(a_1,\dots,a_n)$ if $v\in M_F^0$;
\item $h_v(a_1\cdots a_n) \leq h_v(a_1)+\cdots+h_v(a_n)$.
\end{enumerate}
\end{lemma}

\begin{cor}\label{hv-cor} If $A^{(1)},\dots,A^{(m)}$ are 
$n\times n$-matrices with entries in $F$, then
\begin{enumerate}
\item $h_v\bigl(\det A^{(1)},\dots,\det A^{(m)}\bigr) \leq 
n\cdot\bigl(h_v\bigl(A^{(1)},\dots,A^{(m)}\bigr)+\log n\bigr)$ if $v\in M_F^\infty$;
\item $h_v\bigl(\det A^{(1)},\dots,\det A^{(m)}\bigr) \leq n\cdot h_v\bigl(A^{(1)},\dots,A^{(m)}\bigr)$ if $v\in M_F^0$.
\end{enumerate}
\end{cor}

This follows from (1) and (3)--(5) in the lemma.
Using Hadamard's inequality 
it is possible to
improve the term $\log n$ in (1) slightly, to $\frac{1}{2}\log n$.

The (global) {\em height}\/ of a finite set $S\subseteq F$ is defined
in terms of the local heights:
$$h(S) := \frac{1}{d}\sum_{v\in M_F} h_v(S).$$
The (global) {\em height}\/ of $f_1,\dots,f_n\in F[X]$ 
is  the
global height of its set of coefficients.
The quantity $h(S)$ does not change if the field $F$ is replaced by another
algebraic number field containing the set $S$.
Hence $h$ gives rise to a height function, also denoted by $h$, on
(finite subsets of) the algebraic closure of $\Q$. The product formula
\eqref{Product} implies that $h(a)=h(1/a)$ for all $a\in F^\times$.

We have $h(S')\leq h(S)$ for all subsets $S'\subseteq S$;
hence $h(a)\leq h(S)$ for $a\in S$.
Suppose that $S\neq\{0\}$, and
let $I$ denote the fractional ideal  generated by $S$.
Then
$$
h(S) = \frac{1}{d}\left(\log N(\frak d)+\sum_{v\in M_F^\infty} h_v(S)\right),$$
where $I=\frak b/\frak d$ is a factorization of $I$ with $\frak b,\frak d$
relatively prime ideals of $R$.
In particular, for $0\neq a\in R$ we get
\begin{equation}\label{Height-Norm}
h(a)=h(1/a)=
\frac{1}{d}\left(\log N(a)+\sum_{v\in M_F^\infty} h_v(1/a)\right).
\end{equation}
Moreover, if $S\subseteq R$, then
$$h(S) = \frac{1}{d}\sum_{v\in M_F^\infty} h_v(S).$$
It follows that in this case $h(S)\leq d\cdot \max\bigl\{h(a):a\in S\bigr\}$. 

\begin{example}
For non-zero and relatively prime integers $r,s\in\Z$
we have
$h(r/s)=\max\bigl\{\log\ \abs{r},\log\ \abs{s}\bigr\}$.
In particular $h(r)=\log\ \abs{r}$ for $0\neq r\in\Z$.
\end{example}

From Lemma~\ref{hv-lemma} 
we get the following rules for computing with $h$:
\begin{flalign}
&h(a)=h(-a) \label{h1} \\
&h(a^k) = \abs{k}h(a) \quad\text{for all $k\in\Z$,} \label{h2} \\
&h(a_1+\cdots+a_n) \leq h(a_1,\dots,a_n)+\log n, \label{h6}\\
&h(a_1\cdots a_n) \leq h(a_1)+\cdots+h(a_n), \label{h3}
\end{flalign}
From Corollary~\ref{hv-cor} we obtain the following bound on the height
of determinants of $n\times n$-matrices 
$A^{(1)},\dots,A^{(m)}\in F^{n\times n}$:
\begin{equation}\label{hdet}
h\bigl(\det A^{(1)},\dots,\det A^{(m)}\bigr) \leq n\cdot\bigl(h\bigl(A^{(1)},\dots,A^{(m)}\bigr)+\log n\bigr).
\end{equation}
The following facts will also be used later on:

\begin{lemma}\label{hv}
For all $a\in F^\times$,
$$\sum_{v\in M_F^0, v_{\frak p}(a)>0} 
\log p\cdot v_{\frak p}(a) \leq d\cdot h(a).$$
Here the sum runs over all $v\in M_F^0$ such that
$v_{\frak p}(a)>0$, with 
$\frak p={\frak p}_v$ denoting the prime
ideal of $R$ corresponding to $v$ and $p$ the unique prime number such
that $v|p$.
\end{lemma}
\begin{proof}
We have (using \eqref{h2})
$$
d\cdot h(a)=d\cdot h(1/a)\geq
\sum_{v\in M_F^0} d_v/e_v\cdot\log p\cdot\max\bigl\{0,v_{\frak p}(a)\bigr\}
            \geq \sum_{v\in M_F^0, v_{\frak p}(a)>0} \log p\cdot v_{\frak p}(a)
$$
as claimed.
\end{proof}

It follows that given a non-zero element $a$ of $R$ there are at most
$d\cdot h(a)/\log 2$ 
many absolute values $v\in M_F^0$ such that $v_{\frak p}(a)>0$, where
${\frak p}={\frak p}_v$. Moreover
$\abs{v_{\frak p}(a)}\leq d\cdot h(a)/\log p$ for all $v\in M_F^0$, 
where $v|p$.

\begin{lemma}\label{C-Lemma}
There exists a constant $C_0$, depending only on $F$, 
with the following property:
Given ideals $I$ and $J$ of $R$ with $I$ properly contained in 
$J$, there exists
$a\in J\setminus I$ of height at most $C_0+\frac{1}{d}\log N(J)$.
\end{lemma}
\begin{proof}
Let $\omega_1,\dots,\omega_d$ be a basis for $R$ as $\Z$-module, and
let $$c=d\max\bigl\{ |\omega_1,\dots,\omega_d|_v : v\in M_F^\infty\bigr\}.$$
By Theorem~0 in \cite{Lang}, V, \S{}11 and its proof, 
for any $M_F$-divisor $\frak c$ we have
\begin{equation}\label{Size-Estimate}
\left(\frac{1}{4c}\right)^d ||\frak c||_F < \# L(\frak c) \leq
\max\bigl\{1, 2^{d+2}||\frak c||_F\bigr\}.
\end{equation}
Now let $$t:=4c\cdot N(J)^{1/d}\max\bigl\{1,8/N(I)^{1/d}\bigr\},$$ and 
let
$\frak d$ be the $M_F$-divisor given by
${\frak d}(v) = 1$ if $v\in M_F^0$ and ${\frak d}(v)=t$ if $v\in M_F^\infty$.
We consider the $M_F$-divisors ${\frak c}={\frak c}_I$ and 
${\frak c}'={\frak c}_J\cdot\frak d$ of size
$||\frak c||_F=1/N(I)$ and $||{\frak c}'||_F=t^d/N(J)$, respectively.
By \eqref{Size-Estimate} we have 
$$\# L({\frak c}')>\left(\frac{t}{4c\cdot N(J)^{1/d}}\right)^d =
\max\bigl\{1,8^d/N(I)\bigr\} \geq \max\bigl\{1,2^{d+2}/N(I)\bigr\} \geq
\# L({\frak c}).$$
Therefore $L({\frak c}')\setminus L({\frak c})\neq\emptyset$, that is, 
there exists $a\in J\setminus I$ with $\abs{a}_v \leq t$ for all
$v\in M_F^\infty$. Since $t\leq 32 c N(J)^{1/d}$, it follows that 
$$h(a)=\frac{1}{d}\sum_{v\in M_F^\infty}
d_v \log^+\abs{a}_v \leq \log^+ t \leq C_0 + \frac{1}{d}\log N(J)$$
where $C_0=\log 32+\log c$ is a constant only depending on the number field
$F$.
\end{proof}

\begin{remarks}\ 

\begin{enumerate}
\item The estimate \eqref{Size-Estimate} can be refined to
$$\# L(\frak c)=B_F ||\frak c||_F + O\bigl(||\frak c||_F^{1-1/d}\bigr)
\qquad\text{as $||\frak c||_F\to\infty$,}$$
where $B_F = \frac{2^{d_1}(2\pi)^{d_2}}{\abs{D(F)}^{1/2}}$. 
Here $d_1$ and $d_2$
denote the number of real and complex places of $F$, respectively,
and $D(F)$ denotes the discriminant of $F$.
(See \cite{Lang}, Chapter V, \S{}2, Theorem~1.) 
The constants hidden in the $O$-notation depend on $F$.
We decided to use the
cruder relation \eqref{Size-Estimate}, in order to make the constant $C_0$
as explicit as possible. 
\item By the proof, the constant 
$C_0 =  \log(32d)+h(\omega_1,\dots,\omega_d)$, where
$\omega_1,\dots,\omega_d$ is a basis for the $\Z$-module $R$, has the 
property stated in Lemma~\ref{C-Lemma}. 
A result of 
\cite{Roy-Thunder} implies that one can find such a $\Z$-basis with
$$h(\omega_1,\dots,\omega_d) \leq d\log d + 3d^2\log 2 +
\log\frac{\abs{D(F)}}{i(F)}.$$
Here $i(F)$ denotes the class index of $F$, i.e., the smallest positive
integer such that each ideal class contains an ideal $I$ of $R$
with $N(I)\leq i(F)$.  If $F$ is totally real (i.e., $d_2=0$), then
this bound can be improved to
$$h(\omega_1,\dots,\omega_d) \leq d\log d + \frac{1}{2}d(3d-1)\log 2 +
\frac{1}{2}\log\abs{D(F)}.$$
\end{enumerate}
\end{remarks}

Given any fractional ideal $I$ of $F$ and a non-zero prime ideal $\frak p$ of 
$R$ there exists an element $b$ of $I^{-1}=(R: I)$ such that 
$v_{\frak p}(b)=-v_{\frak p}(I)$. In Section~\ref{Effective-Flatness-Section} we will need the
existence of such $b$ having small height, for integral $I$:

\begin{cor}\label{C-Lemma-Cor}
There exists a constant $C_1$, depending only on $F$, with the following
property: Given an ideal $I=(a_1,\dots,a_n)$ of $R$ and a prime ideal
$\frak p\neq 0$  of $R$ there
exists an element $b$ of $I^{-1}$ 
such that $v_{\frak p}(b)=-v_{\frak p}(I)$ and
$$h(b) \leq C_1\bigl(h(a_1,\dots,a_n)+1\bigr).$$
\end{cor}
\begin{proof}
Let $C_0>0$ be the constant from Lemma~\ref{C-Lemma}, and put $C_1:=C_0+2$.
Let $i\in\{1,\dots,n\}$ be such that $v_{\frak p}(a_i)=v_{\frak p}(I)$, and put
$J:=(a_i)\cdot I^{-1}$, an ideal of $R$. By Lemma~\ref{C-Lemma} there
exists $a\in J\setminus J\cdot \frak p$ with $h(a)\leq C_0+\frac{1}{d}\log N(J)$,
and by \eqref{Height-Norm} 
$$\frac{1}{d}\log N(J) \leq \frac{1}{d}\log N(a_i) \leq h(a_i).$$
Hence the element $b=a/a_i\in I^{-1}$ satisfies
$$h(b) \leq h(a)+h(a_i)\leq C_0+\frac{1}{d}\log N(J) + h(a_i) \leq C_0+2h(a_i)$$
and has the required properties.
\end{proof}

\begin{remarks}\ 

\begin{enumerate}
\item If the number field $F$ is explicitly given, say in terms of its
multiplication table for a $\Z$-basis $\omega_1,\dots,\omega_d$ of $R$,
and the generators $a_1,\dots,a_n$ are also explicitly given (in terms of
their coefficients in the basis $\omega_1,\dots,\omega_d$), then $b\in I^{-1}$
with $v_{\frak p}(b)=-v_{\frak p}(I)$ can be found effectively: By
\cite{Cohen}, p.~202--205 we can compute a basis $b_1,\dots,b_m$
for the $R$-module $I^{-1}$; then $b=b_i$, where $v_{\frak p}(b_i)$ is minimal,
has the required property.
\item For $F=\Q$ the constant $C_1=1$ has the property claimed in the
corollary: Given integers $a_1,\dots,a_n\in\Z$ and a prime number $p$ let 
$b=p^{-\mu}$, where $\mu=\min_i v_p(a_i)$. Then 
$ba_1,\dots,ba_n\in\Z$ and $h(b)=\log p^\mu \leq h(a_1,\dots,a_n)$.
\end{enumerate}
\end{remarks}

\section{Rings of Restricted Power Series}\label{PowerSeriesSection}

Let $\cal O$ be a discrete valuation ring (DVR) with maximal ideal 
$\m=t\cal O$.
We write $v_{\m}\colon {\cal O}\setminus\{0\}\to\N$ for the $\m$-adic
valuation associated to $\cal O$ (normalized so that $v_{\m}(t)=1$). 
We always consider $v_{\m}$ extended to a map $v_{\m}\colon {\cal O}\to\N_\infty$
by $v_{\m}(0):=\infty$,
where $\N_\infty=\N\cup\{\infty\}$ with the usual conventions
$\N<\infty$ and $r+\infty=\infty+r=\infty$ for all $r\in\N_\infty$.
The residue field of
${\cal O}$ is denoted by
$\bar{\cal O}={\cal O}/\m$, with residue homomorphism 
$a\mapsto\bar{a}\colon {\cal O}\to\bar{\cal O}$.

{\it From now until further notice we assume that $\cal O$ is
complete in the $\m$-adic topology on $\cal O$.}\/  
The completion of the polynomial ring ${\cal O}[X]=
{\cal O}[X_1,\dots,X_N]$ with 
respect to the $\m{\cal O}[X]$-adic topology on ${\cal O}[X]$
will be denoted by $\OX = {\cal O}\langle X_1,\dots,X_N \rangle$. 
It may be regarded as a subring
of the ring ${\cal O}[[X]]$ of formal power series over $\cal O$, and is 
called the ring of {\bf restricted power series} with coefficients in $\cal O$.
Its elements are the power series 
$$f=\sum_\nu a_\nu X^\nu\in {\cal O}[[X]] \qquad \text{($a_\nu\in \cal O$ for all $\nu$)}$$
such that $a_\nu\to 0$ (in the $\m$-adic topology on ${\cal O}$)
as $\abs{\nu}\to\infty$. 
Here $\nu=(\nu_1,\dots,\nu_N)$
ranges over all multi-indices in $\N^N$, and $\abs{\nu}=\nu_1+\cdots+\nu_N$. 

The $\m$-adic valuation $v_{\m}\colon {\cal O}\to\N_\infty$
extends to $\OX$
by setting
$$v_{\m}(f)=\min_\nu v_{\m}(a_\nu) \qquad\text{for $f=\sum_\nu a_\nu
X^\nu\in
\OX$.}$$
The map $v_{\m}\colon\OX\to\N_\infty$ is
a valuation on the domain $\OX$, that is, for all $f,g\in\OX$ we have
$v_{\m}(fg)=v_{\m}(f)+v_{\m}(g)$ and
$v_{\m}(f+g)\geq\min\bigl\{v_{\m}(f),v_{\m}(g)\bigr\}$.
(See \cite{BGR}, p.~44, Corollary~2.)
We denote the image of $f\in\OX$ under the canonical surjection 
$\OX\to\OX/t\OX\cong\bar{\cal O}[X]$ by $\bar{f}$.

Suppose from now on that $N\geq 1$, and let $X':=(X_1,\dots,X_{N-1})$.
Canonically $\OXprime\subseteq\OX$, and every element $f\in\OX$
can be written uniquely as 
\begin{equation}\label{f}
f=\sum_{i=0}^\infty f_iX_N^i
\quad\text{with $f_i(X')\in\OXprime$ for all $i\in\N$,}
\end{equation}
where the infinite sum
converges with respect to the $\m\OX$-adic topology on $\OX$. An element 
$f$ of $\OX$, expressed as in \eqref{f}, 
is called {\bf regular in $X_N$ of degree $s\in\N$}
if its reduction $\bar{f}\in \bar{\cal O}[X]$ is unit-monic 
of degree $s$ in $X_N$, that is,
\begin{enumerate}
\item $\bar{f_s}\neq 0$, and
\item $v_{\m}(f_i)>0$ for all $i>s$.
\end{enumerate}
If $f\in\OXprime[X_N]$ is monic of $X_N$-degree $s$
(so that in particular $f$ is regular in $X_N$ of degree $s$, as an element of $\OX$),
then $f$ is called a {\bf Weierstra\ss{} polynomial in $X_N$ of degree $s$.}
For a proof of the following standard facts see, e.g., \cite{BGR}.

%

\begin{lemma}\label{NoetherZpX}
Let $e>1$ and 
suppose that the image of $f\in\OX$ in $\bar{\cal O}[X]$ is non-zero of 
degree $<e$.
Let $T_e\colon\OX\to\OX$ be the $\cal O$-automorphism defined by
\begin{align*}
X_i &\mapsto X_i+X_N^{e^{N-i}} \qquad\text{\rom{(}for $1\leq i<N$\rom{)}}\\
X_N &\mapsto X_N.
\end{align*}
Then $T_e(f)$ is regular in $X_N$ of degree $<e^N$.
\end{lemma}

The ring of restricted power series  has the following fundamental
property:

\begin{theorem}\rom{(Weierstra\ss{} Division Theorem for $\OX$.)}
Let $g\in\OX$ be regular in $X_N$ of degree $s$. Then for each $f\in\OX$
there are uniquely determined elements $q\in \OX$ and $r\in\OXprime[X_N]$
with $\deg_{X_N} r<s$ such that $f=qg+r$.
\end{theorem}
 
In particular, we get
$$
\OX/(g) \iso \OXprime \oplus \OXprime \bar{X_N}\oplus\cdots\oplus\OXprime \bar{X_N}^{s-1}
$$
as $\OXprime$-algebras. (Here, $\bar{X_N}=X_N\bmod g$.)
Applying Weierstra\ss{} Division with $f=X_N^s$, we obtain the
important corollary:

\begin{cor}\label{WP}
\rom{(Weierstra\ss{} Preparation Theorem for $\OX$.)}
Let $g\in\OX$ be regular in $X_N$ of degree $s$. There are a unique
Weierstra\ss{} polynomial $w\in\OXprime[X_N]$ of degree $s$ and a unique unit
$u\in\OX$ such that $g=u\cdot w$.
\end{cor}

From Weierstra\ss{} Preparation it follows that
the ring $\OX$ is Noetherian. 
Here is another useful consequence:

\begin{cor}\label{DivisionByWP}
Let $w\in\OXprime[X_N]$ be a Weierstra\ss{} polynomial. Then the 
inclusion map $\OXprime[X_N]\subseteq \OX$ induces an isomorphism
$$\OXprime[X_N]/w\OXprime[X_N] \overset{\iso}\longrightarrow \OX/w\OX.$$
\end{cor}
\begin{proof}
The surjectivity of the map follows from the existence part of
Weierstra\ss{} Division. For injectivity, we 
have to show: if $fw=g\in\OXprime[X_N]$ for some $f\in\OX$, then 
$f\in\OXprime[X_N]$. This follows by Euclidean Division of $g$ by
the monic polynomial $w$ in $\OXprime[X_N]$, and by the uniqueness statement
in the Weierstra\ss{} Division Theorem.
\end{proof}

Let now $\cal O$ be an arbitrary DVR, not necessarily complete, 
with maximal ideal generated by $t$, and let
$\hat{\cal O}$ be the completion of $\cal O$ in the $\m$-adic topology.
We let $F=\Frac({\cal O})$  be the
fraction field of $\cal O$.
The following lemma and its corollary below will become important in
later sections.

\begin{lemma}\label{SolCompletion}
If a \rom{(}finite\rom{)} system of linear equations over ${\cal O}[X]$ has a
solution in $F[X]$ and in $\OhatX$, then it has a 
solution in ${\cal O}[X]$. 
\end{lemma}
\begin{proof}
For simplicity, we just treat the case of a single linear equation
\begin{equation}\label{SingleEqu}
f_0=f_1y_1+\cdots+f_ny_n \qquad (f_0,f_1,\dots,f_n\in {\cal O}[X]).
\end{equation}
The general case is similar. From a solution 
in $F[X]$ we obtain, after clearing denominators, an integer 
$e\geq 1$ and polynomials
$g_1,\dots,g_n\in {\cal O}[X]$ such that
\begin{equation}\label{SingleEqu1}
t^ef_0 = f_1g_1+\cdots+f_ng_n.
\end{equation}
Now $\OhatX$ is faithfully flat over its subring
$(S_e)^{-1}{\cal O}[X]$, where $S_e$ is the multiplicative set $1+t^e{\cal O}[X]$. 
(See \cite{Greco-Salmon}, Theorems 4.9, 5.1.)
So if \eqref{SingleEqu} is solvable in $\OhatX$, then
there exist $h,h_1,\dots,h_n\in {\cal O}[X]$ with
\begin{equation}\label{SingleEqu2}
(1+t^eh)f_0 = f_1h_1+\cdots+f_nh_n.
\end{equation}
Multiplying \eqref{SingleEqu1} on both sides by $h$ and subtracting
from \eqref{SingleEqu2}, we obtain
$$f_0=f_1(h_1-hg_1)+\cdots+f_n(h_n-hg_n)$$
with $h_1-hg_1,\dots,h_n-hg_n\in {\cal O}[X]$ as desired.
\end{proof}

\begin{cor}\label{SolCompletionCor}
Let $A$ be an $m\times n$-matrix over ${\cal O}[X]$. If 
$$y^{(1)},\dots,y^{(L)}\in \bigl({\cal O}[X]\bigr)^n$$
generate the $F[X]$-module $\Sol_{F[X]}(A)$
of solutions of the homogeneous system 
of linear equations $Ay=0$
in $F[X]$, and
$$z^{(1)},\dots,z^{(M)}\in \bigl({\cal O}[X]\bigr)^n$$
generate the $\OhatX$-module $\Sol_{\OhatX}(A)$ of solutions of
$Ay=0$ in $\OhatX$, then
$$y^{(1)},\dots,y^{(L)},z^{(1)},\dots,z^{(M)}$$
generate the ${\cal O}[X]$-module $\Sol_{{\cal O}[X]}(A)$
of solutions of $Ay=0$ in ${\cal O}[X]$. \qed
\end{cor}

\section{Hermann's Method}\label{HermannSection}

In this section, we first give a presentation of Hermann's method for
constructing generators for the solutions of systems of
homogeneous
linear equations over polynomial rings. We begin by adapting 
this approach so that
it applies to systems of linear equations over any integral domain $D$.
In the next section we will use a variant of Hermann's method in the case where
$D=\OX$ for a complete DVR $\cal O$. 
Here we present the case (treated by Hermann) where $D$ is a polynomial ring 
over a field and
deduce bounds on the degrees of generators for 
syzygy modules. Finally we show how
this method can be modified to solve inhomogeneous systems.

\subsection*{Hermann's method in a general setting}
Let $D$ be an integral domain with fraction field $K$. (Typically, $D$
is a ring of polynomials over an integral domain.)
We consider a homogeneous system of linear equations
\begin{equation}
\tag{I}\label{I}
\begin{bmatrix}
a_{11} & \cdots & a_{1n} \\
\vdots & \ddots & \vdots \\
a_{m1} & \cdots & a_{mn}
\end{bmatrix}
\begin{bmatrix}
y_1 \\
\vdots \\
y_n
\end{bmatrix}
=
\begin{bmatrix}
0 \\
\vdots \\
0
\end{bmatrix}
\end{equation}
with coefficient matrix $A=
(a_{ij})_{\substack{1\leq i\leq m\\
1\leq j\leq n}}$ 
having entries $a_{ij}\in D$.
We are interested in effectively 
finding a set of generators for the module of syzygies $\Sol(A)=\Sol_D(A)$
of $A$. Of course, for this we may assume
$A\neq 0$.
We shall indicate here a reduction of this problem to a similar problem
over a coefficient ring (a quotient of $D$) that is in many cases simpler
than the domain $D$.

Let $r=\rank_K(A)\geq 1$ be the rank of $A$ (considered as a matrix over $K$)
and let $\Delta$ be an $r\times r$-submatrix of $A$ with $\delta=\det\Delta
\neq 0$. After rearranging the order of the equations and permuting the
unknowns $y_1,\dots,y_n$ in \eqref{I} we may assume that $\Delta$ is
the upper left corner of $A$, i.e., $\Delta=(a_{ij})_{1\leq i,j\leq r}$.
Each row $a_i=(a_{i1},\dots,a_{in})$ with $r<i\leq m$ is a $K$-linear
combination
of the first $r$ rows $a_1,\dots,a_r$, so 
\eqref{I} has the same solutions in $D^n$ as
the system:
\begin{equation}
\tag{II}\label{II}
\begin{bmatrix}
a_{11} & \cdots & a_{1n} \\
\vdots & \ddots & \vdots \\
a_{r1} & \cdots & a_{rn}
\end{bmatrix}
\begin{bmatrix}
y_1 \\
\vdots \\
y_n
\end{bmatrix}
=
\begin{bmatrix}
0 \\
\vdots \\
0
\end{bmatrix}
\end{equation}
Changing the notation, we let $r=m$ and $A=(a_{ij})_{\substack{1\leq i\leq r\\
1\leq j\leq n}}$.
So \eqref{II} can now be written as $Ay=0$.
Multiplying both sides of $Ay=0$ on the left by the adjoint of $\Delta$,
\eqref{II} turns into the system
\begin{equation}
\tag{S}\label{S}
\begin{bmatrix}
\delta &        &        &        & c_{1,r+1} & \cdots & c_{1,n} \\
       & \delta &        &        & c_{2,r+1} & \cdots & c_{2,n} \\
       &        & \ddots &        & \vdots    & \ddots & \vdots \\
       &        &        & \delta & c_{r,r+1} & \cdots & c_{rn}
\end{bmatrix}
\begin{bmatrix}
y_1 \\ y_2 \\
\vdots \\
y_n
\end{bmatrix}
=
\begin{bmatrix}
0 \\ 0 \\
\vdots \\
0
\end{bmatrix}
\end{equation}
with $c_{ij},d_i\in D$ for $1\leq i\leq r<j\leq n$, which has the
same solutions in $D^n$ as \eqref{II}, and as \eqref{I}.
We note the following $n-r$
linearly independent solutions of \eqref{S}:
\begin{equation}\label{specialsols}v^{(1)}=\begin{bmatrix}
-c_{1,r+1} \\
\vdots \\
-c_{r,r+1} \\
\delta \\
0 \\
\vdots \\
0
\end{bmatrix},
v^{(2)}=\begin{bmatrix}
-c_{1,r+2} \\
\vdots \\
-c_{r,r+2} \\
0\\
\delta \\
\vdots \\
0
\end{bmatrix}, \dots,
v^{(n-r)}=\begin{bmatrix}
-c_{1,n} \\
\vdots \\
-c_{r,n} \\
0 \\
\vdots \\
0 \\
\delta
\end{bmatrix}
\end{equation}
If $\delta$ is a unit, these vectors form in fact a basis for
$\Sol(A)$. Suppose $\delta$ is not a unit, so $\overline D=D/\delta D\neq 0$.
Then, reducing the coefficients in \eqref{S} modulo $\delta$, the system
\eqref{S} turns into the system
\begin{equation}
\tag{$\overline{\text{S}}$}\label{Sbar}
\begin{bmatrix}
\overline{c_{1,r+1}} & \cdots & \overline{c_{1n}} \\
\vdots & \ddots & \vdots \\
\overline{c_{r,r+1}} & \cdots & \overline{c_{rn}}
\end{bmatrix}
\begin{bmatrix}
y_{r+1} \\
\vdots \\
y_n
\end{bmatrix}
=
\begin{bmatrix}
{0} \\
\vdots \\
{0}
\end{bmatrix}
\end{equation} 
over $\bar{D}$. (Here $\overline{a}$ denotes the image
of $a\in D$ in $\overline{D}$.) 

\begin{lemma}\label{obvious}
Let $z^{(1)},\dots,z^{(M)}\in D^{n-r}$ be such that 
$\bar{z^{(1)}},\dots,\bar{z^{(M)}}\in \bar{D}^{n-r}$ generate the
$\bar{D}$-module of solutions to \eqref{Sbar}.
The vectors $z^{(1)},\dots,z^{(M)}$ may be extended uniquely to
vectors $y^{(1)},\dots,y^{(M)}$ in $D^{n}$
which, together with the solutions of \eqref{I} in \eqref{specialsols}, 
generate $\Sol(A)$.
\end{lemma}

This fact is rather obvious, but what makes it useful is that under
favorable circumstances $\overline D$ is ``simpler'' than $D$. 
(Note however that 
it may happen that $\overline D$ is not a domain anymore.)
Let us consider an example
where this can be exploited. 

\subsection*{Hermann's method for $F[X]$}
Assume that $D$ is a polynomial ring over a field $F$, that is, $D=F[X]=
F[X_1,\dots,X_N]$. Let $N>0$.
Suppose first that $F$ is {\it infinite.}\/ 
In this case, after a linear change of variables, we may assume that
\begin{equation}\label{form}
\delta=uX_N^e + \text{terms of lower $X_N$-degree,}\quad\text{with $e=\deg\delta>0$,
$u\in F^\times$.}
\end{equation}
Then by Euclidean Division each element $\bar{a}\in\bar{D}$ can be uniquely written as
$$\bar{a}=a_0+a_1\bar{X_N}+a_2\bar{X_N}^2+\cdots +a_{e-1}\bar{X_N}^{e-1}$$
with $a_0,\dots,a_{e-1}\in F[X']=F[X_1,\dots,X_{N-1}]$. 
In particular, each coefficient
$\bar{c_{ij}}$ in \eqref{Sbar} can be written in this way.
Note that $\deg_{X'}a_i\leq\deg_X a$ for all $0\leq i<e$.
Let us also
write each unknown $y_j$  in \eqref{Sbar}, for $r<j\leq n$, as
$$y_j=y_{j0}+y_{j1}\bar{X_N}+\cdots + y_{j,e-1}\bar{X_N}^{e-1}$$
with new unknowns $y_{jk}$ ($r<j\leq n$, $0\leq k<e$) ranging over $D'=F[X']$. Each product
$\bar{c_{ij}}y_j$ in \eqref{Sbar} can then be written as
$$\beta_0(y_{j0},\dots,y_{j,e-1})+\beta_1(y_{j0},\dots,y_{j,e-1})\bar{X_N}+\cdots+
\beta_{e-1}(y_{j0},\dots,y_{j,e-1})\bar{X_N}^{e-1},$$
where each $\beta_k$ is a linear form in $y_{j0},\dots,y_{j,e-1}$ with
coefficients in $D'$. From this, it is routine to construct a homogeneous 
system
of $r(e-1)$ linear equations in the $e(n-r)$ unknowns $y_{jk}$ over $D'$
whose solutions in $D'$ are in one-to-one correspondence with the
solutions of $\eqref{Sbar}$ in $\bar{D}$. 

\subsection*{Computing degree bounds}
For the sake of obtaining ``good'' bounds on the degrees of solutions, we modify the general construction sketched 
above, exploiting some more special features of $F[X]$. 
Put $d=\deg_{X_N} A$.
Write each $a_{ij}$ as
\begin{equation}\label{aij}
a_{ij}=a_{ij0}+a_{ij1}X_N+\cdots+a_{ijd}X_N^d
\end{equation}
with $a_{ijk}\in F[X']$, and also each unknown $y_j$ as
\begin{equation}\label{yj}
y_j=y_{j0}+y_{j1}X_N+\cdots+y_{j,rd-1}X_N^{rd-1}
\end{equation}
with new unknowns $y_{jk}$ ranging over $F[X']$. Then
the $i$-th equation in \eqref{II} yields $(r+1)d$ equations
$$\sum_{l=0}^k \sum_{j=1}^n a_{ijl} y_{j,k-l} = 0,\qquad\qquad
0\leq k<(r+1)d,$$
where we put $a_{ijl}:=0$ for $l>d$ and $y_{i,l}:=0$ for $l\geq rd$ .
In this way, we obtain a new system
\begin{equation}\tag{I$'$}\label{Iprime}
A'y'=0,
\end{equation}
where $A'$ is an $\bigl(rd(r+1)\bigr)\times (nrd)$-matrix with entries in $D'$
and
\begin{equation}\label{yprime}
y'=\bigl[
y_{1,0}, \dots, y_{1,rd-1}, \dots, y_{n,0}, \dots, y_{n,rd-1}\bigr]^{\tr},
\end{equation}
whose solutions in $D'$ are in one-to-one 
correspondence with the solutions of \eqref{II} in $D$
of $X_N$-degree $<rd$. 
Note that the entries of $A'$  are still of degree (in $X'$)
at most $\deg_{X} A$. If $N>1$, then we can repeat the
same procedure with \eqref{Iprime} instead of \eqref{I}, etcetera, until 
we obtain a (huge) homogeneous system of linear
equations over $F$. We can (effectively) find a finite set of
generators for the $F$-vector space of solutions to this system, and
reversing the process above, we
obtain a finite set of generators for the original system
\eqref{I}: Suppose we have already found
a finite set of generators
for the $D'$-submodule $\Sol_{D'}(A')$ of $(D')^{nrd}$,
where $A'$ is the matrix constructed from $A$ as above.
That is, we have finitely many solutions $y^{(1)},\dots,y^{(M')}$ of
\eqref{I} such that each solution to \eqref{I} of $X_N$-degree $<rd$
is a linear combination of $y^{(1)},\dots,y^{(M')}$. 
The solutions in \eqref{specialsols} together with
$y^{(1)},\dots,y^{(M')}$ form a set of generators for $\Sol(A)=\Sol_D(A)$: 
Given
any solution $y=[y_1,\dots,y_n]^{\tr}\in\Sol(A)$ we can divide
each $y_j$, $j=n-r+1,\dots,n$ by $\delta$:
$$y_j=Q_{j-r}\delta+R_{j-r} \qquad (j=n-r+1,\dots,n)$$
with $Q_1,\dots,Q_{n-r}\in F[X]$ and $R_1,\dots,R_{n-r}\in F[X]$ of 
$X_N$-degree $<e$. Then
$$z=y-Q_1v^{(1)}-\cdots-Q_{n-r}v^{(n-r)}=\bigl[h_1, \dots, h_r, R_1, \dots, 
R_{n-r}\bigr]^{\tr}$$
is also a solution to \eqref{S}, with $h_1,\dots,h_r\in F[X]$. Now
$$\delta h_i = -(c_{i,r+1}R_1+\cdots+c_{in}R_{n-r}) \qquad
\text{for $i=1,\dots,r$,}$$
where the right-hand sides have $X_N$-degree $<rd+e$. Hence $\deg_{X_N} h<
rd$ and therefore $\deg_{X_N} z<rd$. It follows that $z$ is a 
$D$-linear combination of 
$y^{(1)},\dots,y^{(M')}$, so $y$ is a $D$-linear combination of
$y^{(1)},\dots,y^{(M')},v^{(1)},\dots,v^{(n-r)}$ as claimed.

Let $\alpha=\alpha(N,d,m)$ be the smallest natural number such that for all
infinite fields $F$, a system of $m$ homogeneous
linear equations \eqref{I} over $D=F[X]=F[X_1,\dots,X_N]$
with all $\deg a_{ij}$ bounded from above by $d$
is generated by the solutions
of degree $\leq\alpha$. (By the considerations above, $\alpha(N,d,m)$ exists.)
The derived system \eqref{Iprime}
consists of at most $dm(m+1)$ equations in at most 
$dn^2$ unknowns, and $\deg_{X'}(A')\leq d$. 
From a set of generators of the solutions
to \eqref{Iprime} of degree $\leq d'$
we can produce a set of generators of the solutions to
\eqref{I} of degree $\leq d'+md$.
We get the relation
$$\alpha(N,d,m) \leq \alpha\bigl(N-1,d,dm(m+1)\bigr)+md$$
for $N>0$. Noting that $\alpha(0,d,m)=0$
for all $d$, $m$, we find that
$$\alpha(N,d,m)\leq (m+1)d+\bigl((m+1)d\bigr)^2 + \cdots +
\bigl((m+1)d\bigr)^{2^{N-1}}\leq (2md)^{2^{N}}.$$
If $F$ is any field, possibly finite, we work over $F'=F(T)$,
an infinite field. Here, $T$ is an indeterminate 
distinct from $X_1,\dots,X_N$. Given $y\in \bigl(F[T,X])^n$ write
$y=y(0)+y(1)T+y(2)T^2+\cdots$ (a finite sum)
with $y(k)\in \bigl(F[X]\bigr)^n$ for all $k$. If
${\cal G}$ is a generating set for $\Sol_{F'[X]}(A)$ consisting of
elements of $\bigl(F[T,X]\bigr)^n$, then the collection of $y(k)$, where
$y\in{\cal G}$ and $k\in\N$, generates $\Sol_{F[X]}(A)$.
To sum up, we have shown the classical result:

\begin{theorem}\rom{(Hermann \cite{Hermann}, Seidenberg \cite{Seidenberg1}.)}\label{HermannMainThm}
For every polynomial ring $D=F[X_1,\dots,X_N]$ over a field $F$ and $A\in
D^{m\times n}$ of degree $\leq d$, the solution module $\Sol_D(A)$ of the 
homogeneous system
$Ay=0$ is generated by the solutions 
of degree at most $\beta(N,d,m)=(2md)^{2^{N}}$. \qed 
\end{theorem}


\subsection*{Hermann's method for inhomogeneous systems}
Let again $D$ be a domain with fraction field $K$.
Given an $m\times n$-matrix $A=(a_{ij})$
with entries $a_{ij}\in D$, we are now interested in
determining for each column vector $b=[b_1,\dots,b_m]^{\tr}\in
D^m$ whether the system
\begin{equation}
\tag{I$_b$}\label{Ib}
\begin{bmatrix}
a_{11} & \cdots & a_{1n} \\
\vdots & \ddots & \vdots \\
a_{m1} & \cdots & a_{mn}
\end{bmatrix}
\begin{bmatrix}
y_1 \\
\vdots \\
y_n
\end{bmatrix}
=
\begin{bmatrix}
b_1 \\
\vdots \\
b_m
\end{bmatrix}
\end{equation}
(or: $Ay=b$) is solvable for some $y=[y_1,\dots,y_n]^{\tr}
\in D^n$, and if
it is, effectively finding such a solution. 
Similarly to the case of homogeneous equations, this problem can be reduced
to an analogous problem over a quotient of $D$:
As above let $\Delta$ be an $r\times r$-submatrix of $A$ with 
$\delta=\det\Delta\neq 0$, where $r=\rank_K(A)\geq 1$. 
Again we may assume that $\Delta=(a_{ij})_{1\leq i,j\leq r}$.
Each row $a_i=(a_{i1},\dots,a_{in})$ with $r<i\leq m$ is a $K$-linear
combination
$$a_i=\sum_{\varrho=1}^r \lambda_{i\varrho} a_\varrho \qquad
(\lambda_{i\varrho}\in K)$$
of the first $r$ rows $a_1,\dots,a_r$. 
So a {\it necessary condition}\/ for \eqref{Ib} to
have a solution in $D^n$ is that
\begin{equation}
\label{NC}\tag{NC}
b_i=\sum_{\varrho=1}^r \lambda_{i\varrho} b_\varrho \qquad
\text{for $r<i\leq m$.}
\end{equation}
(That is, $\rank_K(A)=\rank_K(A,b)$.)
Assume \eqref{NC} holds. Then \eqref{Ib} has the same solutions in $D^n$ as
the system:
\begin{equation}
\tag{II$_b$}\label{IIb}
\begin{bmatrix}
a_{11} & \cdots & a_{1n} \\
\vdots & \ddots & \vdots \\
a_{r1} & \cdots & a_{rn}
\end{bmatrix}
\begin{bmatrix}
y_1 \\
\vdots \\
y_n
\end{bmatrix}
=
\begin{bmatrix}
b_1 \\
\vdots \\
b_r
\end{bmatrix}
\end{equation}
Changing the notation, we let $r=m$, 
so \eqref{IIb} can now be written as $Ay=b$.
Multiplying both sides of $Ay=b$ on the left by the adjoint $\Delta^{\operatorname{ad}}$
of $\Delta$,
\eqref{IIb} turns into the system
\begin{equation}
\tag{S$_b$}\label{Sb}
\begin{bmatrix}
\delta &        &        &        & c_{1,r+1} & \cdots & c_{1,n} \\
       & \delta &        &        & c_{2,r+1} & \cdots & c_{2,n} \\
       &        & \ddots &        & \vdots    & \ddots & \vdots \\
       &        &        & \delta & c_{r,r+1} & \cdots & c_{rn}
\end{bmatrix}
\begin{bmatrix}
y_1 \\ y_2 \\
\vdots \\
y_n
\end{bmatrix}
=
\begin{bmatrix}
d_1 \\ d_2 \\
\vdots \\
d_r
\end{bmatrix}
\end{equation}
(with $c_{ij},d_i\in D$ for $1\leq i\leq r<j\leq n$) which has the
same solutions in $D^n$ as \eqref{IIb}, and as \eqref{Ib}. 
Clearly, a sufficient condition for
\eqref{Sb} to have a solution $y=[y_1,\dots,y_n]^{\tr}\in D^n$ 
is that $d_1,\dots,d_r$ are
each divisible by $\delta$. This will be the case 
if $\delta$ is a unit. A solution to \eqref{Sb} (and
hence to \eqref{Ib}) is then given by 
$$y_j=\begin{cases} d_j/\delta &\text{for $1\leq j\leq r$,}\\
                    0 &\text{for $d<j\leq n$.}\end{cases}$$ 
Suppose $\delta$ is not a unit, so $\bar{D}=D/\delta D\neq 0$.
Then, reducing the
coefficients in \eqref{Sb} modulo $\delta$, the system \eqref{Sb} turns into
\begin{equation}
\tag{$\overline{\text{S}_b}$}\label{Sbarb}
\begin{bmatrix}
\overline{c_{1,r+1}} & \cdots & \overline{c_{1n}} \\
\vdots & \ddots & \vdots \\
\overline{c_{r,r+1}} & \cdots & \overline{c_{rn}}
\end{bmatrix}
\begin{bmatrix}
y_{r+1} \\
\vdots \\
y_n
\end{bmatrix}
=
\begin{bmatrix}
\overline{d_1} \\
\vdots \\
\overline{d_r}
\end{bmatrix}
\end{equation}
over $\overline{D}$. The key fact here is the following (similar to
Lemma~\ref{obvious}):

\begin{lemma}
Any $[y_{r+1},\dots,y_n]^{\tr}\in D^{n-r}$ with the property that 
$\bigl[\overline{y_{r+1}},\dots,\overline{y_n}\bigr]^{\tr}$ is 
a solution of the reduced system \eqref{Sbarb} can be
augmented uniquely to a solution $$y=\bigl[y_1,\dots,y_r,y_{r+1},\dots,y_n
\bigr]^{\tr}\in D^n$$
of \eqref{Sb}, and hence of \eqref{Ib}.
\rom{(In particular, \eqref{Ib} is solvable in $D$ if and only if \eqref{Sbarb} 
is solvable in $\bar D$.)}
\end{lemma}

In the case where $D=F[X]$ is a polynomial ring over a field $F$, 
we can again modify this reduction somewhat to facilitate
the computation of bounds. Suppose that $N>0$
and $F$ is infinite. Then, after applying a linear change of variables,
we may assume that $\delta$ has the form \eqref{form}. 
By Euclidean Division we write each $b_i$ 
as
$$b_i=\delta f_i+g_i \qquad \text{with $f_i,g_i\in D$, $\deg_{X_N} g_i<e$.}$$
The solutions of \eqref{IIb} in $D^n$ are in one-to-one
correspondence with the solutions in $D^n$ of the system
\begin{equation}
\tag{III$_b$}\label{IIIb}
\begin{bmatrix}
a_{11} & \cdots & a_{1n} \\
\vdots & \ddots & \vdots \\
a_{r1} & \cdots & a_{rn}
\end{bmatrix}
\begin{bmatrix}
y_1 \\
\vdots \\
y_n
\end{bmatrix}
=
\begin{bmatrix}
g_1 \\
\vdots \\
g_r
\end{bmatrix}
\end{equation}
with the same coefficient matrix $A$ as \eqref{IIb}. 
To see this, let $$f=\left[\begin{smallmatrix}f_1\\ \vdots\\ f_r
\end{smallmatrix}\right], \quad g=\left[\begin{smallmatrix}g_1\\ \vdots\\ g_r
\end{smallmatrix}\right],\quad \text{and} \quad
h=\left[\begin{smallmatrix} \Delta^{\operatorname{ad}}f \\ 0 \\ \vdots \\ 0
\end{smallmatrix}\right]\in D^n.$$ 
Note that $b=\delta f+g$ and $Ah=\delta f$; so
$y\in D^n$ is a solution to
\eqref{IIb} if and only if $y-h\in D^n$ is a solution to \eqref{IIIb}. 
Moreover, if all $a_{ij}$ and $b_i$ have degree $\leq d$, and if
\eqref{IIb} is solvable in $D^n$, then \eqref{IIIb} even has a solution
in $D^n$ of $X_N$-degree $<rd$. In order
to prove this, suppose $y=[y_1,\dots,y_n]^{\tr}\in D^n$ is a solution to 
\eqref{IIb}. The polynomial 
$\delta$, each $c_{ij}$ and each $d_i$ have degree at most $rd$.
Subtracting from $y$
appropriate multiples of the solutions $v^{(1)},\dots,v^{(n-r)}$ 
(see \eqref{specialsols}) to the homogeneous system $Ay=0$
associated with \eqref{IIb}, if necessary, we may assume that
$\deg_{X_N} y_j<e\leq rd$ for $j=r+1,\dots,n$.  
Multiplying the equation $A(y-h)=g$ on both sides from 
the left by the adjoint 
$\Delta^{\operatorname{ad}}$ of $\Delta$, we get, for $j=1,\dots,r$:
$$\delta(y_j-h_j)=\sum_{k=r+1}^n y_kc_{jk}+
\text{(terms of $X_N$-degree $<e+rd$)}.$$
It follows that $\deg_{X_N}
\bigl(\delta(y_j-h_j)\bigr)<e+rd$ and thus $\deg_{X_N}
(y_j-h_j)<rd$. So $y-h$ is a solution to \eqref{IIIb} 
of $X_N$-degree $<rd$ as required.

Write each $g_i$ as
$$g_i=g_{i0}+g_{i1}X_N+\cdots+g_{i,e-1}X_N^{e-1}$$
with $g_{i0},\dots,g_{i,e-1}\in F[X']$, each $a_{ij}$ 
in the form \eqref{aij}, and also each unknown $y_j$ as in \eqref{yj}.
Comparing the
coefficients of equal powers of $X_N$ on both sides,
the $i$-th equation in \eqref{IIIb} yields $(r+1)d$ equations
$$\sum_{l=0}^k \sum_{j=1}^n a_{ijl} y_{j,k-l} = g_{ik},\qquad\qquad
0\leq k<(r+1)d,$$
with  $a_{ijl}:=0$ for $l>d$, $y_{i,l}:=0$ for $l\geq rd$,
$g_{il}:=0$ for $l\geq e$.
We get a new system
\begin{equation}\tag{I$_b'$}\label{Ibprime}
A'y'=b',
\end{equation}
where $A'$ is an $\bigl(rd(r+1)\bigr)\times (nrd)$-matrix with entries in $D'$,
$b'$ is an $\bigl(rd(r+1)\bigr)$-column vector with components
from $D'$, and $y'$ as in \eqref{yprime},
whose solutions in $D'$ are in one-to-one 
correspondence with the solutions of \eqref{IIIb} in $D$
of $X_N$-degree $<rd$.
So starting with a system \eqref{I} over $D=F[X_1,\dots,X_N]$ we 
have constructed a system
\eqref{Ibprime} over $D'=F[X_1,\dots,X_{N-1}]$ which is, assuming \eqref{NC}, in some sense equivalent to it. 
Note that $\deg_{X'}(A',b')\leq d$.

Associated to \eqref{Ibprime} we have the necessary condition
\begin{equation}\label{NCprime}\tag{NC$'$}
\rank_{K'}(A')=\rank_{K'}(A',b') \qquad \text{(where $K'=\Frac(D')$)}
\end{equation}
for its solvability in $D'$. So if $N>1$ and \eqref{NCprime} holds, 
then we can repeat the procedure with \eqref{Ibprime}, 
until we obtain a system of linear
equations over $K$. We can (effectively) decide whether this system has
a solution over $K$, and if it does,  find one,
e.g., by Gaussian Elimination.  
Eventually we
obtain a solution $y\in D^n$ of the original system
\eqref{Ib} with $\deg y \leq \beta(N,d,m)=(2md)^{2^{N}}$, where $d=\deg(A,b)$.

If $F$ is a finite field, we again  work over the infinite field $F'=F(T)$.
The algorithm described above allows to test whether the system \eqref{Ib}
has a solution $y'=[y_1',\dots,y_n']^{\tr}\in\bigl(F'[X]\bigr)^n$, 
and if it is, effectively obtain such a solution with
$\deg_X y' \leq \beta(N,d,m)$.
Since the coefficients of the $y_j'$ solve a certain
system of linear equations involving the coefficients of the $a_{ij}$ and
the $b_i$, we
can also find a solution $y$ {\it in $F[X]$}\/
with the $\deg y$ majorized by the same bound. 
This shows:

\begin{theorem}\rom{(Hermann \cite{Hermann}, Seidenberg \cite{Seidenberg1}.)}
\label{HermannMainThm2}
For every polynomial ring $D=F[X_1,\dots,X_N]$ over a field $F$ 
and $A\in D^{m\times n}$, $b\in D^m$
of degree $\leq d$, if the system of linear equations
$Ay=b$ has a solution in $D^n$ 
then it has such a solution of degree at most $(2md)^{2^{N}}$. \hfill\qed
\end{theorem}

The following is a consequence of the preceding
theorem and Cramer's Rule:

\begin{cor}\label{Cor-Hermann}
Let $R$ be a domain with fraction field $F=\Frac(R)$, and let
$D=R[X]=R[X_1,\dots,X_N]$.
Given a finitely generated
submodule $M$ of the free $D$-module $D^m$,
there exists a non-zero $\delta\in R$
with the property that 
\begin{equation}\label{delta}
v\in MF[X] \quad\Longleftrightarrow\quad
\delta v\in M \qquad\text{for all $v\in D^m$.}
\end{equation} 
If $R$ is computable, then $\delta$ can be computed elementary
recursively \rom{(}in the ring operations of $R$\rom{)} from given
generators for $M$. \qed
\end{cor}

\begin{remark}
Theorems~\ref{HermannMainThm} and \ref{HermannMainThm2} above remain true
for a polynomial ring $D=R[X_1,\dots,X_N]$ over a von Neumann regular ring
$R$. This follows easily from
the fact that for any von Neumann regular ring $R$ there exists a
faithfully flat embedding $R\to S$ into a direct product $S$ of fields.
\end{remark}

\section{Effective Flatness}\label{Effective-Flatness-Section}

The purpose of this section is to prove Theorem~B from
the introduction, in a more general setting.
A ring $R$ is called {\bf hereditary} if every ideal 
of $R$ is projective (as an $R$-module). 
A domain $R$ is hereditary if and only if $R$ is
a Dedekind domain. (\cite{Glaz}, p.~27.) 
A domain $R$ is called {\bf almost Dedekind}
if every localization $R_{\frak m}$ of $R$ at a maximal ideal $\frak m$
of $R$ is a DVR. (See \cite{Gilmer}, p.~434.) Somewhat more generally,
we shall call a ring $R$ {\bf almost hereditary} if
the ring of fractions $\Frac(R)$ of $R$ is von Neumann
regular and $R_{\frak m}$ is a DVR for every maximal ideal $\frak m$ of
$R$. If $R$ is hereditary, then $R$ is almost hereditary. 
(\cite{Glaz}, pp.~27--28.)
There exist
examples of domains which are almost Dedekind but not Dedekind; 
see \cite{Gilmer}, pp.~516--518. With this terminology, we have:

\begin{theorem}\label{EffectiveFlatnessTheorem}
Let $R$ be an almost hereditary ring and $A=(a_{ij})\in D^{m\times n}$,
$A\neq 0$, where $D=R[X_1,\dots,X_N]$. 
The module of solutions to $Ay=0$ in $D$ is generated by elements of
degree $\leq (2m\deg A)^{2((N+1)^N-1)}$.
\end{theorem}

(Since an almost hereditary ring is semihereditary and hence coherent, see
\cite{Glaz}, p.~128,
finitely many such generators will suffice.)

As a first step in the proof of this theorem, we show an easy
local-global result:

\begin{lemma}\label{First-Lemma}
Let $R$ be a ring 
with ring of fractions $F=\Frac(R)$,
and let $M$ be an $R[X]$-submodule of $R[X]^n$.
For each maximal ideal $\frak m$ of $R$ let
$v_{\frak m}^{(1)},\dots,v_{\frak m}^{(K_{\frak m})}\in M$
generate the $R_{\frak m}[X]$-submodule $MR_{\frak m}[X]$ of $R_{\frak m}[X]^n$
generated by \rom{(}the canonical image of\rom{)} $M$. Then
$v_{\frak m}^{(1)},\dots,v_{\frak m}^{(K_{\frak m})}$,
where $\frak m$ ranges over all maximal ideals of $R$,
generate the $R[X]$-module $M$.
\end{lemma}
\begin{proof}
Let $y\in M$. Then 
for any maximal ideal $\frak m$ of $R$ there
exists $\delta_{\frak m}
\in R\setminus\frak m$ and $b_{1,\frak m},\dots, b_{K_{\frak m},\frak m}
\in R[X]$ such that
\begin{equation}\label{y2}
\delta_{\frak m} y = b_{{\frak m},1}v_{\frak m}^{(1)}+\cdots+b_{{\frak m},K_{\frak m}}v_{\frak m}^{(K_{\frak m})}.
\end{equation}
The various
$\delta_{\frak m}$, where $\frak m$ ranges over all maximal ideals of $R$,
generate the unit ideal of $R$. Hence there exist
maximal ideals ${\frak m}_1,\dots,{\frak m}_k$ of $R$ (for some $k\in\N$)
and $c_1,\dots,c_k\in R$ such that
$$1=c_1 \delta_{{\frak m}_1} + \cdots + c_k \delta_{{\frak m}_k}.$$
Therefore
$$y=c_1 (\delta_{{\frak m}_1}y) + \cdots + c_k 
(\delta_{{\frak m}_k}y).$$
Together with \eqref{y2} this shows that $y$
is an $R[X]$-linear combination of the
$v_{\frak m}^{(j)}$.
\end{proof}

\begin{remark}
Suppose that $u^{(1)},\dots,u^{(K)}\in M$ generate the $F[X]$-module
$MF[X]$, and let $\delta\in (M':M)\cap R$, where
$M'$ is the $R[X]$-submodule of $M$ generated by $u^{(1)},\dots,u^{(K)}$.
Similarly to the proof of the lemma one shows that
$u^{(1)},\dots,u^{(K)}$ together with
$v_{\frak m}^{(1)},\dots,v_{\frak m}^{(K_{\frak m})}$, where $\frak m$
ranges over all maximal ideals of $R$ containing $\delta$, suffice to generate
$M$.
\end{remark}

Let now $R$ be an almost hereditary ring
and $0\neq A=(a_{ij})\in D^{m\times n}$, where $D=R[X]$, $X=(X_1,\dots,X_N)$.
Then $F=\Frac(R)$ is von Neumann regular, and $R_{\frak m}$ is a DVR, for
every maximal ideal $\frak m$ of $R$.
By virtue of the lemma applied to $M=\Sol_{D}(A)$, it suffices to
find 
$$v_{\frak m}^{(1)},\dots,v_{\frak m}^{(K_{\frak m})}\in\Sol_D(A)$$ 
generating $MR_{\frak m}[X]=\Sol_{R_{\frak m}[X]}(A)$, for each maximal ideal $\frak m$ of $R$,
with $v_{\frak m}^{(j)}$ 
of ``small'' degree. 
For the construction of the $v_{\frak m}^{(j)}$ 
we may use Corollary~\ref{SolCompletionCor},
since $R_{\frak m}$ is a DVR. Hence, 
given a maximal ideal $\frak m$ of $R$ we need to find 
\begin{enumerate}
\item $y_{\frak m}^{(1)},\dots,y_{\frak m}^{(L_{\frak m})}\in\Sol_{R_{\frak m}[X]}(A)$ generating $\Sol_{\Frac(R_{\frak m})[X]}(A)$ and
\item $z_{\frak m}^{(1)},\dots,z_{\frak m}^{(M_{\frak m})}\in\Sol_{R_{\frak m}[X]}(A)$ generating $\Sol_{\hat{R_{\frak m}}\langle X\rangle}(A)$,
\end{enumerate}
with $y_{\frak m}^{(i)}$ and $z_{\frak m}^{(j)}$ 
of degree $\leq (2m\deg A)^{2((N+1)^N-1)}$. 
By Hermann's Theorem \ref{HermannMainThm} from the last section we obtain
$y_{\frak m}^{(i)}$ satisfying (1), of degree bounded by
$$(2m\deg A)^{2^N}\leq (2m\deg A)^{2((N+1)^N-1)} \qquad \text{(for $N>0$).}$$
 The existence of the
$z_{\frak m}^{(j)}$ is a consequence of the following
{\it effective flatness}\/ result applied to the DVR ${\cal O}=R_{\frak m}$:

\begin{prop}\label{EffectiveFlatness}
Let $\cal O$ be a DVR with maximal ideal $\m$ and
$\m$-adic completion $\hat{\cal O}$, and $A=(a_{ij})\in
\bigl({\cal O}[X]\bigr)^{m\times n}$, $A\neq 0$.
There exist solutions $z^{(1)},\dots,z^{(M)}\in\Sol_{{\cal O}[X]}(A)$ of degree
at most $(2m\deg A)^{2((N+1)^N-1)}$ which generate $\Sol_{\OhatX}(A)$.
\end{prop}

In proving this proposition we proceed by induction on $N$, 
following Hermann's method as in
the proof of Theorem \ref{HermannMainThm}, with $F[X]$ replaced by
$\OhatX$ and Weierstra\ss{} Division for $\OhatX$ in place of Euclidean
Division for $F[X]$. However this procedure breaks down if $\delta \bmod
t=0$ for all $r\times r$-minors $\delta$ of $A$, since then Weierstra\ss{}
Division by $\delta$ is {\em inapplicable.}\/ To overcome this
obstacle, we shall first transform our system 
\begin{equation}
\tag{I}\label{Iflat}
\begin{bmatrix}
a_{11} & \cdots & a_{1n} \\
\vdots & \ddots & \vdots \\
a_{m1} & \cdots & a_{mn}
\end{bmatrix}
\begin{bmatrix}
y_1 \\
\vdots \\
y_n
\end{bmatrix}
=
\begin{bmatrix}
0 \\
\vdots \\
0
\end{bmatrix}
\end{equation} 
into an equivalent system for which $\delta\bmod t\neq 0$ for a suitable
$r\times r$-minor $\delta$ of the new coefficient matrix. For this, by 
removing superfluous rows from $A$ we may of course 
assume that the rows of $A$ are linearly independent over the fraction 
field $F(X)$ of ${\cal O}[X]$, i.e., $m=r=\rank_{F(X)}(A)\geq 1$.
Let $\Delta$ be an $r\times r$-submatrix of $A$ such that
$v_{\m}(\det\Delta)$ is {\em minimal}\/ among all $r\times r$-submatrices
of $A$. Without loss of generality, $\Delta=(a_{ij})_{1\leq i,j\leq
r}$. As in Section~\ref{HermannSection}, consider now the system
\begin{equation}\label{Sflat}
\tag{S}
\begin{bmatrix}
\delta &        &        &        & c_{1,r+1} & \cdots & c_{1,n} \\
       & \delta &        &        & c_{2,r+1} & \cdots & c_{2,n} \\
       &        & \ddots &        & \vdots    & \ddots & \vdots \\
       &        &        & \delta & c_{r,r+1} & \cdots & c_{rn}
\end{bmatrix}
\begin{bmatrix}
y_1 \\ y_2 \\
\vdots \\
y_n
\end{bmatrix}
=
\begin{bmatrix}
0 \\ 0 \\
\vdots \\
0
\end{bmatrix}
\end{equation}
which is obtained by multiplying both sides of \eqref{Iflat} from the
left with the adjoint of $\Delta$. It has the same solutions as
\eqref{Iflat} in any domain extending ${\cal O}[X]$. Here, $\delta=\det\Delta$,
and the $c_{ij}$ are certain signed $r\times r$-minors of $A$. 
In particular, $v_{\m}(c_{ij})\geq v_{\m}(\delta)$ for all $i$, $j$, by
choice of $\Delta$. 
We have the $n-r$
linearly independent solutions
\begin{equation}\label{Special}
v^{(1)}=\begin{bmatrix}
-c_{1,r+1} \\
\vdots \\
-c_{r,r+1} \\
\delta \\
0 \\
\vdots \\
0
\end{bmatrix},
v^{(2)}=\begin{bmatrix}
-c_{1,r+2} \\
\vdots \\
-c_{r,r+2} \\
0\\
\delta \\
\vdots \\
0
\end{bmatrix}, \dots,
v^{(n-r)}=\begin{bmatrix}
-c_{1,n} \\
\vdots \\
-c_{r,n} \\
0 \\
\vdots \\
0 \\
\delta
\end{bmatrix}
\end{equation}
to the homogeneous system \eqref{Sflat}. Put  $\mu=v_{\m}(\delta)$
and $u^{(k)}=t^{-\mu}v^{(k)}\in \bigl({\cal O}[X]\bigr)^n$ for $k=1,\dots,n-r$.
If $N=0$, then $t^{-\mu}\delta$ is a unit in $\cal O$,
so the solutions $u^{(1)},\dots,u^{(n-r)}$
form a basis of $\Sol_{\cal O}(A)$ and hence of $\Sol_{\hat{\cal O}}(A)$
(since $\hat{\cal O}$ is flat over $\cal O$).
Suppose now that $N>0$. 
We let 
$e=r\deg A+1$ and put $b_{ij}=T_e(a_{ij})$, where $T_e$
is the $\hat{\cal O}$-automorphism of $\OhatX$ 
defined in Lemma~\ref{NoetherZpX}. 
Then the system $By=0$, where $B=(b_{ij})
\in \bigl({\cal O}[X]\bigr)^{m\times n}$, has the same rank $r$ as \eqref{Iflat}, and 
$y\in\OhatX^n$ is a solution to \eqref{I} if and only if $T_e(y)$ is a
solution to $By=0$.
Dividing all coefficients $\delta$ and $c_{ij}$ 
in \eqref{S} by $t^\mu$ and applying $T_e$ to the resulting system,
we obtain a system 
\begin{equation}\label{Sflate}\tag{S$_e$}
\begin{bmatrix}
\varepsilon &        &        &        & d_{1,r+1} & \cdots & d_{1,n} \\
       & \varepsilon &        &        & d_{2,r+1} & \cdots & d_{2,n} \\
       &        & \ddots &        & \vdots    & \ddots & \vdots \\
       &        &        & \varepsilon & d_{r,r+1} & \cdots & d_{rn}
\end{bmatrix}
\begin{bmatrix}
y_1 \\ y_2 \\
\vdots \\
y_n
\end{bmatrix}
=
\begin{bmatrix}
0 \\ 0 \\
\vdots \\
0
\end{bmatrix}
\end{equation}
which has the same solutions, in any domain extending ${\cal O}[X]$, as 
$By=0$, where $d_{ij}\in {\cal O}[X]$ for all $i,j$ and
$\varepsilon\in {\cal O}[X]$ is regular in $X_N$ of some degree $s<e^N$.
This system has the $n-r$ linearly independent solutions
$w^{(1)},\dots,w^{(n-r)}$, where $w^{(k)}=T_e\bigl(u^{(k)}\bigr)$ 
for $k=1,\dots,n-r$.
Let $d\deg_{X_B} B$, so
$d<e^N$. (Note that $\deg_{X'}(b_{ij})\leq\deg_{X'}(a_{ij})$ for all $i,j$.)
Write
$$b_{ij}=b_{ij0}+b_{ij1}X_N+\cdots+b_{ijd}X_N^d$$
with $b_{ij0},\dots,b_{ijd}\in {\cal O}[X']$, 
and each unknown $y_j$ as
$$y_j=y_{j0}+y_{j1}X_N+\cdots+y_{j,rd-1}X_N^{rd-1}$$
with new unknowns $y_{jk}$ ($1\leq j\leq n$, $0\leq k<rd$) ranging over
$\OhatXprime$.  The $i$-th equation in $By=0$ may then be written as
$$\sum_{l=0}^k \sum_{j=1}^n b_{ijl}y_{j,k-l}=0,\qquad\qquad
0\leq k<(r+1)d,$$
where we put $b_{ijl}:=0$ for $l>d$ and $y_{i,l}:=0$ for $l\geq rd$.
This gives rise to a system over ${\cal O}[X']$:
\begin{equation}\label{AprimeF}
A'y'=0,
\end{equation}
consisting of $rd(r+1)$ homogeneous equations in the $nrd$ 
unknowns $y'=(y_{jk})$,
whose solutions in $\OhatXprime$ are in one-to-one correspondence with
the solutions 
$y\in\bigl(\OhatXprime[X_N]\bigr)^n$ to $By=0$
with $\deg_{X_N} y<rd$.
From a finite set of generators of $\Sol_{{\cal O}[X']}(A')$ we thus obtain
finitely many column vectors
$$y^{(1)},\dots,y^{(M')} \in \bigl({\cal O}[X]\bigr)^n$$
with the following property: each $y^{(i)}$ is a solution to
``$By=0$'' of $X_N$-degree $<rd$, 
and each solution $y
\in\bigl(\OhatXprime[X_N]\bigr)^n$ to this system of linear
equations with $\deg_{X_N} y<rd$ is an
$\OhatXprime$-linear combination of $y^{(1)},\dots,y^{(M')}$.
Consider now the solutions
\begin{equation}\label{SolsA}
u^{(1)},\dots,u^{(n-r)},T_e^{-1}\bigl(y^{(1)}\bigr),\dots,
T_e^{-1}\bigl(y^{(M')}\bigr)\in\bigl({\cal O}[X]\bigr)^n
\end{equation}
to \eqref{I}. We show:

\begin{lemma}\label{SolsALemma}
The vectors in \eqref{SolsA} generate the $\OhatX$-module $\Sol_{\OhatX}(A)$.
\end{lemma}
\begin{proof}
Suppose that $x\in\bigl(\OhatX\bigr)^n$
is any solution to $Ay=0$, and let $y=T_e(x)$, a solution to $By=0$.
Since $\varepsilon$ is regular in $X_N$ of
degree $s$, we can write, by Weierstra\ss{} Division in $\OhatX$:
$$y_j=Q_{j-r} \varepsilon + R_{j-r}, \qquad (j=n-r+1,\dots,n)$$
with $Q_1,\dots,Q_{n-r}\in\OhatX$ and $R_1,\dots,R_{n-r}\in\OhatXprime[X_N]$
of $X_N$-degree $<s$. Then
$$z=y-Q_{1}w^{(1)}-\cdots-Q_{n-r}w^{(n-r)}=\bigl[h_1,\dots,h_r,
R_1,\dots,R_{n-r}\bigr]^{\tr}$$
is also a solution to \eqref{Sflate}, with $h_1,\dots,h_r\in\OhatX$. Let
$U\in\OhatX$ be a unit and $W\in\OhatXprime[X_N]$ be a Weierstra\ss{}
polynomial such that $\varepsilon=UW$. Since $\varepsilon$ 
is polynomial in $X_N$,
by Lemma~\ref{DivisionByWP} we also have $U\in\OhatXprime[X_N]$.
The degree of $\varepsilon$ in $X_N$ is $\leq rd$, and the degree of $W$ in
$X_N$ is $s$; hence $U$ is of degree $\leq rd-s$ in $X_N$. Moreover,
\begin{equation}\label{wuhi}
W(Uh_i)=\varepsilon h_i = -(d_{i,r+1}R_{1}+\cdots+d_{in}R_{n-r}) \in
\OhatXprime[X_N]
\end{equation}
for $i=1,\dots,r$.
Since $W$ is monic in $X_N$, it follows that
$Uh_i\in\OhatXprime[X_N]$. Put $z=Uz'\in\bigl(\OhatXprime[X_N]\bigr)^n$, a
solution to \eqref{Sflate}. We claim that all entries of $z$
have $X_N$-degree $<rd$: To see this note that
$\deg_{X_N} R_i < s$ for $i=1,\dots,n-r$
and $\deg_{X_N} U\leq d-s$; hence the last $n-r$ entries $UR_1,\dots,UR_{n-r}$
of $z$ are of $X_N$-degree $<rd$. For the first $r$ entries $Uh_1,\dots,Uh_r$
use that the right-hand side of \eqref{wuhi} has $X_N$-degree $<rd+s$;
since $\deg_{X_N} W=s$ we get $\deg_{X_N} Uh_i<rd$ for all $i=1,\dots,r$.
It follows that $z'$, and hence $z$, is an $\OhatX$-linear combination of
$y^{(1)},\dots,y^{(K')}$. Since $U$ is a unit in $\OhatX$, the solution
$y$ can be expressed as an $\OhatX$-linear combination of the column
vectors
$$w^{(1)},\dots,w^{(n-r)},y^{(1)},\dots,y^{(M')}\in\bigl({\cal O}[X]\bigr)^n.$$
Hence the solution $x=T_e^{-1}(y)$ 
to our original equation 
\eqref{I} is an $\OhatX$-linear combination of the vectors in \eqref{SolsA}
as claimed.
\end{proof}
\begin{remark}
We can bound the degrees of the solutions in \eqref{SolsA}:
We have $\deg u^{(k)}\leq r\deg A$ for $k=1,\dots,n-r$ and
$\deg_{X'} T_e^{-1}\bigl(y^{(i)}\bigr)\leq\deg_{X'} y^{(i)}$ for 
$i=1,\dots,M'$. Moreover $\deg_{X_N} y^{(l)}<rd< re^N$ and thus 
$$\deg_{X_N} T_e^{-1}\bigl(y^{(i)}\bigr) \leq e^{N-1}\deg y^{(i)}
\leq e^{N-1}\bigl(\deg_{X'} y^{(i)}+re^N\bigr)
$$
for $i=1,\dots,M'$.
\end{remark}
Starting with \eqref{I} we successively obtain 
equivalent homogeneous matrix equations
\begin{align}
\tag{H$_N$}\label{HN} A^{(N)}y^{(N)} &= 0 \\
\notag\vdots& \\
\tag{H$_\nu$}\label{Hnu} A^{(\nu)} y^{(\nu)} &= 0 \\
\notag\vdots& \\
\tag{H$_0$}\label{H0} A^{(0)} y^{(0)} &= 0,
\end{align}
where $0\leq\nu\leq N$, $A^{(\nu)}$ is an $m(\nu)\times n(\nu)$-matrix with entries in the polynomial ring
${\cal O}[X_1,\dots,X_\nu]$ and $$y^{(\nu)}=\left[y^{(\nu)}_{1},\dots, y^{(\nu)}_{n(\nu)}\right]^{\tr}$$
is a vector of unknowns ranging over $\hat{\cal O}\langle
X_1,\dots,X_\nu\rangle$. 
So the {\em initial}\/ equation \eqref{HN} is just
$Ay=0$, and if $\nu>0$, then the system
(H$_{\nu-1}$) 
is obtained from \eqref{Hnu} by the
procedure described above (passage from $A$ to $A'$). We have
$$m(\nu)\leq m(\nu+1)\bigl(m(\nu+1)+1\bigr)e(\nu+1)^{\nu}$$
for all $\nu=0,\dots,N-1$, where
$e(\nu)=m(\nu)\deg A^{(\nu)}+1$. It follows that
$$e(\nu) \leq m(\nu+1)\bigl(m(\nu+1)+1\bigr)e(\nu+1)^\nu \deg A^{(\nu)} + 1$$
Using that $\deg A^{(\nu)}\leq \deg A$ we get the estimate
\begin{equation}\label{enu}
e(\nu) \leq (m\deg A+1)^{(N+1)^{N-\nu}}
\end{equation}
for all $\nu=0,\dots,N$. Let ${\cal B}(0)\subseteq {\cal O}^{n(0)}$ be a 
finite system of
generators of $\Sol_{\cal O}(A^{(0)})$, and for every $\nu=1,\dots,N$ let  
${\cal B}(\nu)\subseteq {\cal O}[X_1,\dots,X_\nu]^{n(\nu)}$ be 
a system of generators for the module of solutions to
(H$_\nu$) in $\hat{\cal O}\langle X_1,\dots,X_\nu\rangle$, with ${\cal B}(\nu)$ 
constructed from ${\cal B}(\nu-1)$ according to the process described
above. For $\nu=0,\dots,N$ let $\gamma(\nu)$ be the maximal degree of an
element of ${\cal B}(\nu)$.
Clearly $\gamma(0)=0$, and by the remark following
Lemma~\ref{SolsALemma} we have
$$\gamma(\nu) \leq e(\nu)^{\nu-1}\bigl(\gamma(\nu-1)+m(\nu)e(\nu)^\nu\bigr)+\gamma(\nu-1).$$
The right-hand side can be further estimated from above by
$$e(\nu)^{\nu-1}\bigl(2\gamma(\nu-1)+m(\nu)e(\nu)^\nu\bigr) \leq e(\nu)^{2\nu-1}\bigl(\gamma(\nu-1)+m(\nu)\bigr).$$
Hence we get $$\gamma(\nu)+1 \leq e(\nu)^{2\nu}\bigl(\gamma(\nu-1)+1\bigr)$$
for all $\nu=1,\dots,N$. It follows that
$$\gamma(N)+1 \leq e(N)^{2N} e(N-1)^{2(N-1)} \cdots e(1)^2,$$
and hence, using \eqref{enu}:
$$\gamma(N) \leq \bigl(m\deg A+1\bigr)^\varrho$$
where $\varrho=2\sum_{i=0}^{N-1} (N+1)^i (N-i)$. It is easy to see that
$\varrho\leq 2\bigl((N+1)^N-1\bigr)$. 
Hence every element of ${\cal B}(N)$ has degree 
$\leq (2m\deg A)^{2((N+1)^N-1)}$, finishing the proof of
Proposition~\ref{EffectiveFlatness}, and thus of
Theorem~\ref{EffectiveFlatnessTheorem}. \qed

\begin{remarkNumbered}\label{Algorithm-Remark}
As a consequence of Theorem~\ref{EffectiveFlatnessTheorem},
if $R$ is an almost Dedekind domain that is 
syzygy-solvable, then there exists an (impractical) algorithm
which, given an $m\times n$-matrix $A$ with entries in $D=R[X]$, constructs
a finite collection of generators for $\Sol_{D}(A)$.
If $R=\Z$, or more generally, a computable principal ideal domain, 
we can also turn the proof of the theorem into such an algorithm:
We first find
generators $u^{(1)},\dots,u^{(K)}\in\Sol_D(A)$ for $\Sol_{F[X]}(A)$, where
$F=\Frac(R)$.
(See the remark following Theorem~\ref{HermannMainThm}.)
By  Corollary~\ref{Cor-Hermann}
we then can compute $0\neq \delta\in R$ such that 
for every solution $y\in\bigl(R[X]\bigr)^n$ of $Ay=0$,
$\delta y$ is an $R[X]$-linear combination of
$u^{(1)},\dots,u^{(K)}$. Hence $\delta\in (M':M)$, where $M=\Sol_D(A)$,
$M'=Du^{(1)}+\cdots+Du^{(K)}$.
For every prime factor $\pi$ of $\delta$ we now follow the 
inductive procedure outlined in the proof of
Proposition~\ref{EffectiveFlatness} to construct
generators $v_\pi^{(1)},\dots,v_\pi^{(K_\pi)}\in\Sol_D(A)$ for
$\Sol_{\widehat{R_{(\pi)}}\langle X\rangle}(A)$. 
By the remark following Lemma~\ref{First-Lemma}, the solutions
$u^{(1)},\dots,u^{(K)}$ together with the 
$v_\pi^{(1)},\dots,v_\pi^{(K_\pi)}$ (with $\pi$ ranging over the prime
factors of $\delta$), generate $\Sol_D(A)$.
\end{remarkNumbered}

Sometimes Theorem~\ref{EffectiveFlatnessTheorem} 
still holds for rings which are not almost hereditary:

\begin{cor}\label{IntegralClosure}
Let $R$ be an integrally closed
almost Dedekind domain, and let $S$ be the integral closure of
$R$ inside an algebraic closure of the fraction field $F$ of $R$. 
Let $A$ be an $m\times n$-matrix with entries in $S[X]=S[X_1,\dots,X_N]$.
Then $\Sol_{S[X]}(A)$ is generated by elements of degree
at most $(2m\deg A)^{2((N+1)^{N}-1)}$.
\end{cor}
\begin{proof}
Let $F'$ be a finite field extension of $F$ containing all the
coefficients of the entries of $A$, and let $R'$ be the integral closure
of $R$ in $F'$. Then $R'$ is almost Dedekind. (See \cite{Gilmer}, (36.1).)
Since $R'$ is a Pr\"ufer domain and $S$ 
a torsion-free $R'$-module, $S$ is flat over $R'$.
The claim now follows from Theorem~\ref{EffectiveFlatnessTheorem}.
\end{proof}

The corollary applies  to $R=\Z$ (so $S=$ the ring of
all algebraic integers).

\subsection*{Application 1: bounds for module-theoretic operations}

Let $R$ be an almost Dedekind domain with fraction field
$F=\Frac(R)$. We can exploit 
Theorem~\ref{EffectiveFlatnessTheorem} to establish
bounds for some basic operations on finitely generated submodules
of free modules over $D=R[X]=R[X_1,\dots,X_N]$. We say that
a finitely generated $D$-submodule of $D^m$ 
is {\bf of type $d$} (where $d\in\N$) if it is generated by
vectors of degree $\leq d$. 

\begin{prop}\label{ModuleOps}
Let $M$ and $M'$ be finitely generated submodules 
of the free $D$-module $D^m$ of type $d$.
Then the $D$-modules $(M'F[X])\cap D^m$ and $M\cap M'$ and the
ideal $(M':M)$ are of type $(2md)^{2^{O(N^2)}}$.
\end{prop}
\begin{proof}
Let $M=Dv^{(1)}+\cdots+Dv^{(n)}$ and $M'=Dw^{(1)}+\cdots+Dw^{(p)}$ with 
$v^{(i)},w^{(j)}\in D^m$ of
degree $\leq d$. Let $0\neq \delta\in R$ satisfy \eqref{delta}.
To find generators for the $D$-module $(M'F[X])\cap D^m$, we first find a 
finite set of generators $z^{(1)},\dots,z^{(K)}\in D^{p+m}$ for the
$D$-module of solutions to the system of homogeneous equations
$$w^{(1)}y_1+\cdots+w^{(p)}y_p+(-\delta e^{(1)}) y_{p+1}+\cdots+ (-\delta e^{(m)})y_{p+m}=0.$$
Here $e^{(1)},\dots,e^{(m)}$ denote the unit vectors in $D^m$.
Then clearly the $K$ vectors consisting of the last $m$ entries of
$z^{(1)},\dots,z^{(K)}$ generate
$(M'F[X])\cap D^m$ and are of type
$(2md)^{2((N+1)^N-1)}$. This shows (1). 
Similarly, in order to find generators for
$M\cap M'$ it suffices to find generators for the $D$-module of solutions to
the system of homogeneous equations
$$v^{(1)}y_1+\cdots+v^{(n)}y_n=w^{(1)}y_{n+1}+\cdots+w^{(p)}y_{n+p}.$$
Moreover we have 
$$(M': M)=(M': Dv^{(1)})\cap\cdots\cap(M': Dv^{(n)}),$$
and if $u^{(1)},\dots,u^{(q)}\in D^m$ generate $M'\cap Dv$, where 
$v=[v_1,\dots,v_m]^{\tr}\in D^m$,
then 
$$(M': Dv) = \bigcap_{j=1}^m \bigl(u_j^{(1)}/v_j,\dots,u_j^{(q)}/v_j\bigr).$$
Here $a/0:=1$ for all $a\in R$. From this parts (2) and (3) follow easily.
\end{proof}

\begin{remarks}\ 

\begin{enumerate}
\item If $R$ is syzygy-solvable, then generators for the $D$-modules
$(M' F[X])\cap D^m$ and $M\cap M'$ and for the ideal $(M': M)$
can be computed elementary recursively (in the basic operations of $R$)
from given generators for $M$ and $M'$.
This follows from the proof of the proposition and 
Remark~\ref{Algorithm-Remark} above.
\item By Corollary~\ref{IntegralClosure}, 
the proposition remains true if $R$ is replaced
by the ring of algebraic integers.
\end{enumerate}
\end{remarks}

\subsection*{Application 2: a criterion for primeness}

The following lemma is well-known; we leave the proof to the reader.

\begin{lemma}\label{PrimeLemma}
Let $R$ be a ring and $I$ be an ideal of $R[X]$. Then $I$ is
prime if and only if the image of $I$ in $(R/I\cap R)[X]$ is prime.
If $R$ is an integral domain with fraction field $F$ and $I\cap R=(0)$, then
$I$ is
prime if and only if $IF[X]$ is prime and $IF[X]\cap R[X]=I$.
\end{lemma}

As a consequence, we obtain a test for primeness of an ideal in $\Z[X]$.
Given $f_1,\dots,f_n\in\Z[X]$ we choose $0\neq\delta=\delta(f_1,\dots,f_n)\in\Z$
satisfying \eqref{delta} for $R=\Z$ and $M=$ the ideal generated by
$f_1,\dots,f_n$. 

\begin{cor}\label{PrimeCor}
An ideal $I=(f_1,\dots,f_n)$ of 
$\Z[X]$ is prime if and only if one of the following holds:
\begin{enumerate}
\item $I\Q[X]$ is prime and $I\Q[X]\cap\Z[X]=I$, or
\item there exists a prime factor $p$ of $\delta(f_1,\dots,f_n)$ 
such that $p\in I$ and the image of $I$ in
$\F_p[X]$ is a prime ideal.
\end{enumerate}
\end{cor}

Combining Corollary~\ref{PrimeCor} with Proposition~\ref{ModuleOps} and
a result from \cite{Schmidt-Goettsch} we get a criterion for the
primeness of an ideal of $\Z[X]$ which is polynomial in the degrees of
the generators:

\begin{prop}\label{PrimeProp}
There exists $\varrho=\varrho(N)\in\N$ such that for each ideal
$I$ of $\Z[X]$ of type $d$, the following
is true: $I$ is prime if and only if $1\notin I$, and for all
$f,g\in\Z[X]$ of degree $\leq d^\varrho$, 
if $fg\in I$, then $f\in I$ or $g\in I$.
\end{prop}
\begin{proof}
By Proposition~\ref{ModuleOps}~(1), for all ideals
$I$ of $\Z[X]$ of type $d$, 
the ideal $I\Q[X]\cap\Z[X]$ of $\Z[X]$ is
of type $\tau=(2d)^{2((N+1)^N-1)}$. Moreover by \cite{Schmidt-Goettsch} 
there exists $\varrho'=\varrho'(N)\in\N$ such
that for each field $F$ and each ideal $J$ of
$F[X]$ of type $d$, we have: $J$ is
prime if and only if $1\notin J$, and for all $f,g\in F[X]$ of degree
$\leq d^{\varrho'}$, if $fg\in J$, then $f\in J$ or $g\in J$. We claim that
$\varrho=\max\bigl\{4((N+1)^N-1),\varrho'\bigr\}$
has the required properties. For this,
let $I=(f_1,\dots,f_n)$ be an ideal of $\Z[X]$ of type $d$, and put
$\delta=\delta(f_1,\dots,f_n)$.
Suppose $1\notin I$ and
$fg\in I \Rightarrow \text{$f\in I$ or $g\in I$,}$ for all
$f,g\in\Z[X]$ of degree $\leq d^\varrho$. Then $I\Q[X]$ is prime, since
$d^\varrho\geq d^{\varrho'}$. 
Let $f\in I\Q[X]\cap\Z[X]$ be of degree
at most $\tau$. Then 
$\delta f\in I$, and hence $f\in I$ or
$p\in I$ for some prime divisor $p$ of $\delta$. 
Suppose $f\in I$ for all $f\in I\Q[X]\cap\Z[X]$ of degree $\leq
\tau$; then $I\Q[X]\cap\Z[X]=I$, and by 
Corollary~\ref{PrimeCor}~(1) it follows that 
$I$ is prime. If, on the other hand, we have $p\in I$ for some prime divisor 
$p$ of $\delta$, then the
image of $I$ in $\F_p[X]$ is a prime ideal. By Corollary~\ref{PrimeCor}~(2) 
it follows again that
$I$ is prime, as desired.
\end{proof}

It is clear that Corollary~\ref{PrimeCor} and Proposition~\ref{PrimeProp}
hold, mutatis mutandis,
for any PID $R$ with fraction field $F$ in place of $\Z$ and $\Q$,
respectively.

\section{Height Bounds}\label{Height-Section}

Throughout this section we let $F$ be a number field
and  $R={\cal O}_{F}$ the ring of integers of $F$. 
Let $A=(a_{ij})$ be an $m\times n$-matrix
with entries $a_{ij}$ in $D=R[X]=R[X_1,\dots,X_N]$. 
Let $d=\deg A$ and $h=h(A)$.
As was shown in the previous section, 
we can explicitly bound the degrees of generators for 
the $D$-module $\Sol_{D}(A)$ in terms of $d$, $m$ and $N$.
We now want to bound the heights of those generators
in a similar fashion (in terms of $d$, $h$, $m$ and $N$).

\subsection*{The local case}
Let $\frak p \neq 0$ be a  prime ideal of $R$, and
${\cal O}:=R_{\frak p}$ (a DVR). We first investigate the height of
generators for $\Sol_{{\cal O}[X]}(A)$ and begin
with the case $N=0$:

\begin{lemma}\label{N0-Lemma}
Suppose that $a_{ij}\in R$ for all $i,j$, and let $r=\rank_F (A)$.
The $\cal O$-module $\Sol_{\cal O}(A)$ of solutions in ${\cal O}^n$ 
to the system of 
homogeneous linear equations $Ay=0$ is generated by $n-r$
many vectors whose height is bounded by 
$$C_2\cdot r(h+\log r+1).$$
Here $C_2$ is a constant only depending on $F$.
\end{lemma}
\begin{proof}
Let $v\in M_F$ denote the place of $F$ associated with $\frak p$,
so ${\frak p}={\frak p}_v$. 
We may assume that 
$\det\Delta \neq 0$, where $\Delta=(a_{ij})_{1\leq i,j\leq r}$
(after permuting the unknowns in our system $Ay=0$ if necessary).
In fact, we may assume that the $\frak p$-adic valuation
$\mu:=v_{\frak p}(\det\Delta)$ of $\det\Delta$ is minimal among all
$r\times r$-submatrices of $A$, 
cf.~the proof of Proposition~\ref{EffectiveFlatness}.
Now $Ay=0$ has the same solutions in any domain extending $\cal O$ as
the system \eqref{S} 
obtained from $Ay=0$ by multiplying both sides from the
left with the adjoint of $\Delta$ (see Section~\ref{HermannSection}). 
The entries $\delta=\det\Delta$ and 
$c_{ij}$ ($1\leq i\leq r<j\leq n$) of the coefficient matrix of \eqref{S}
are certain signed $r\times r$-minors of $A$. Let
$v^{(1)},\dots,v^{(n-r)}$ are the $n-r$ be the linearly
independent solutions to $Ay=0$ listed in \eqref{specialsols}. 
By \eqref{hdet} we have $h\bigl(v^{(k)}\bigr)\leq r(h+\log r)$
for $k=1,\dots,n-r$.
Corollary~\ref{C-Lemma-Cor} implies that there exists an element $b$ of $F$ 
such that $v_{\frak p}(b)=-\mu$,  $b v^{(k)}\in R^n$ for all $k=1,\dots,n-r$, 
and $h(b) \leq C_1r(h+\log r+1)$. Here $C_1>0$ is a constant
which only depends on $F$.
The vectors $bv^{(1)},\dots,bv^{(n-r)}\in R^n$ generate
$\Sol_{\cal O}(A)$ and are bounded in height by
$C_2r(h+\log r+1)$, with $C_2=2C_1$.
\end{proof}

\begin{remark}
Note that the vectors $v^{(1)},\dots,v^{(n-r)}\in\Sol_{R}(A)$ as in the
proof of the lemma generate
$\Sol_F(A)$ and satisfy $h(v^{(k)})\leq r(h+\log r)$
for $k=1,\dots,n-r$.
Moreover, the element $0\neq \delta\in R$ of height
$h(\delta) \leq r(h+\log r)$
has the property that $\delta\in (M': M)$, where $M=\Sol_{R}(A)$ 
and $M'=Rv^{(1)}+\cdots+Rv^{(n-r)}$.
\end{remark}

We now consider the general case $N \geq 0$.
By Proposition~\ref{EffectiveFlatness}, the solution module 
$\Sol_{{\cal O}[X]}(A)$ is generated by solutions
$y=\left[y_1,\dots, y_n\right]^{\tr}$ with
$y_1,\dots,y_n\in {\cal O}[X]$ of degree at most $\gamma
=\gamma(N,d,m):=(2md)^{2((N+1)^N-1)}$.
Write $$y_j=\sum_{\abs{\nu}\leq \gamma} y_{j,\nu}X^\nu$$
with variables $y_{j,\nu}$ ranging over $\cal O$  and
$$a_{ij}=\sum_{\abs{\mu}\leq d} a_{ij,\mu}X^\mu$$
with $a_{ij,\mu}\in\cal O$,  where $1\leq i\leq m$, $1\leq j\leq n$
and $\nu,\mu\in\N^N$, $\abs{\nu} \leq \gamma$, $\abs{\mu}\leq d$.
A polynomial in $X_1,\dots,X_N$ of degree 
at most $d$ has at most $M(N,d)=\binom{N+d}{N}$ monomials. Hence the solutions
(in ${\cal O}[X]$) of every equation
$$a_{i1}y_1+\cdots+a_{in}y_n=0 \qquad (1\leq i\leq m)$$
are in one-to-one correspondence with the solutions  (in $\cal O$)
of the system consisting of the $M(N,\gamma+d)$ homogeneous equations 
$$\sum_j \sum_{\mu+\nu=\lambda} a_{ij,\mu}y_{j,\nu} = 0 \qquad (\abs{\lambda}
\leq\gamma+d)$$
in the $n\cdot M(N,\gamma)$ many variables $y_{j,\nu}$, with coefficients
in $\cal O$. So the entire system $Ay=0$, 
with coefficients in ${\cal O}[X]$, may be replaced by a certain 
homogeneous system of $m\cdot M(N,\gamma+d)$ equations in the variables 
$y_{j,\nu}$, having coefficients in $\cal O$. Applying the lemma above to 
the new system and using the estimate
$$
m\cdot M(N,\gamma+d) \leq m\cdot (\gamma+d+1)^N =
\bigl( 2m(d+1)\bigr)^{2^{O(N^2+1)}},
$$
we get the following result, with $C_2$ as above. 

\begin{prop}\label{HeightProp}
For any $N\geq 0$,
the ${\cal O}[X]$-module $\Sol_{{\cal O}[X]}(A)$ is generated by solutions of
degree at most
$(2md)^{2^{O(N^2)}}$
and height at most
\begin{equation}\label{HeightPropEq}
C_2\cdot\bigl(2m(d+1)\bigr)^{2^{O(N^2+1)}}(h+1).
\end{equation}
Here $C_2$ is a constant only depending on $F$. \qed
\end{prop}

With $\beta=\beta(N,m,d)=(2md)^{2^{N}}$ we have
$$
m\cdot M(N,\beta+d) \leq m\cdot (\beta+d+1)^N =
\bigl( 2m(d+1)\bigr)^{2^{O(N+1)}}.
$$
Using this estimate as well as
Theorem~\ref{HermannMainThm} (in place of Proposition~\ref{EffectiveFlatness})
and the remark following Lemma~\ref{N0-Lemma}, one obtains
a result similar to Proposition~\ref{HeightProp}:

\begin{lemma}\label{HeightLemma}
The $F[X]$-module $\Sol_{F[X]}(A)$ is generated by vectors
$u^{(1)},\dots,u^{(K)}\in\Sol_{{\cal O}[X]}(A)$ of
degree at most
$(2md)^{2^N}$
and height at most
\begin{equation}\label{HeightLemmaEq}
\bigl(2m(d+1)\bigr)^{2^{O(N+1)}}(h+1).
\end{equation}
Moreover, 
there exists $0\neq\delta\in (M':M)$ of height bounded by 
\eqref{HeightLemmaEq},
where $M=\Sol_{D}(A)$ and $M'=Du^{(1)}+\cdots+Du^{(K)}$.
\qed
\end{lemma}

\subsection*{The global case}
Proposition~\ref{HeightProp} and Lemma~\ref{HeightLemma} now 
imply the existence of generators of $\Sol_{D}(A)$ of small height:

\begin{cor}\label{HeightCor}
The $D$-module $\Sol_{D}(A)$ can be generated by  solutions of
degree at most
$(2md)^{2^{O(N^2)}}$
and of height at most
$$C_2\cdot\bigl(2m(d+1)\bigr)^{2^{O(N^2+1)}}(h+1).$$
Here $C_2$ is a constant only depending on $F$.
\end{cor}
\begin{proof}
By  Lemma~\ref{HeightLemma} we find  
$u^{(1)},\dots,u^{(K)}\in\Sol_{D}(A)$
with the following properties:
$u^{(1)},\dots,u^{(K)}$  generate the $F[X]$-module $\Sol_{F[X]}(A)$, 
$\deg u^{(k)}\leq (2md)^{2^N}$ for all $k$, and
$h(u^{(k)})$ is bounded by \eqref{HeightLemmaEq}, for each $k$. 
Moreover we find an element $0\neq \delta\in (M':M)$ 
of height bounded by \eqref{HeightLemmaEq}, where
$M=\Sol_{D}(A)$ and $M'=Du^{(1)}+\cdots+Du^{(K)}$.
For every maximal ideal $\frak m$ of $R$ we find generators
$v_{\frak m}^{(1)},\dots,v_{\frak m}^{(K_{\frak m})}\in\Sol_D(A)$
of $\Sol_{R_{\frak m}[X]}(A)$
having degree at most $(2md)^{2^{O(N^2)}}$ and
height bounded by \eqref{HeightPropEq}.
By the remark  following
Lemma~\ref{First-Lemma}, the vectors $u^{(1)},\dots,u^{(K)}$ and
$v_{\frak m}^{(1)},\dots,v_{\frak m}^{(K_{\frak m})}$, where $\frak m$ ranges
over all maximal ideals of $R$ containing $\delta$, generate $\Sol_{D}(A)$.
\end{proof}

\begin{remark}
The number of generators of $\Sol_{D}(A)$
can be bounded in a similar way:
If $\delta$ is a unit in $R$, then $u^{(1)},\dots,u^{(K)}$ generate $\Sol_D(A)$,
and $K\leq n\cdot M(N,\beta+d) = n\bigl(2m(d+1)\bigr)^{2^{O(N+1)}}$.
In general, 
by the remark after
Lemma~\ref{hv}, there are at most $[F:\Q]\cdot h(\delta)/\log 2$ many
maximal ideals of $R$ containing $\delta$.
So we have at most 
$$
n\cdot M(N,\gamma+d)\cdot \bigl(1+ [F:\Q]\cdot h(\delta)/\log 2\bigr) = 
n\cdot [F:\Q] \cdot
\bigl(2m(d+1)\bigr)^{2^{O(N^2+1)}} (h+1)
$$
generators in total.
\end{remark}

\section{Ideal Membership}\label{Ideal-Section}

In this section we use the results obtained so far to
give a proof of Theorem~A from the introduction. We 
begin by studying ideal membership problems of a special form.

\subsection*{B\'ezout identities}
Let $R$ be a ring, $f_1,\dots,f_n\in R[X]$, and $d=\max_i\deg f_i$.
We call a representation of $1$ as a linear combination
\begin{equation}\label{Bezout}
1=f_1g_1+\cdots+f_ng_n
\end{equation}
of $f_1,\dots,f_n$ with coefficients $g_1,\dots,g_n\in R[X]$ a 
{\bf B\'ezout identity} for $f_1,\dots,f_n$ in $R[X]$.
If $R=F$ is a field, then from
Hermann's Theorem~\ref{HermannMainThm2} it follows that
$1\in (f_1,\dots,f_n)F[X]$ if and only
if there exist $g_1,\dots,g_n\in F[X]$ of degree $\leq (2d)^{2^N}$
satisfying the B\'ezout identity \eqref{Bezout}. 
By the effective version of Hilbert's Nullstellensatz due to 
Koll\'ar \cite{Kollar}, this bound may be improved substantially:
if $1\in (f_1,\dots,f_n)F[X]$, then
there are $g_1,\dots,g_n\in F[X]$ of degrees $\leq (3d)^N$
satisfying \eqref{Bezout}. For $F=\Q$ this
means: if $1\in (f_1,\dots,f_n)\Q[X]$, then there are
$\delta\in\Z\setminus\{0\}$ and $g_1,\dots,g_n\in\Z[X]$ of degree
$\leq (3d)^N$ with
$$\delta=f_1g_1+\cdots+f_ng_n.$$
We have the following bound for the size of
$\delta$, obtained
along the lines of Lemma~\ref{HeightLemma}
(i.e., Cramer's rule).
From now on, $F$ denotes a number field.

\begin{lemma}\label{delta-lemma}
Suppose that $R={\cal O}_F$ is the ring of integers of $F$.
If we have $1\in (f_1,\dots,f_n)F[X]$, then
$$\delta = f_1g_1+\cdots+f_ng_n$$
for some $g_1,\dots,g_n\in R[X]$ of degree $\leq (3d)^N$ and some
$\delta\in R$, $\delta\neq 0$, of height at most
$$\bigl(2(d+1)\bigr)^{O(N^2+1)}\bigl(h(f_1,\dots,f_n)+1\bigr).$$
\end{lemma}

\begin{remark}
In fact, the height of the denominator $\delta$ can be bounded
in terms of  $N$, $d$, $n$ and the height of $f_1,\dots,f_n$
by a bound which is single-exponential in $d$ and linear in
$h(f_1,\dots,f_n)$, while at the
same time retaining a single-exponential bound on the degrees of the $g_j$.
See, e.g., \cite{Krick-Pardo-2} or \cite{Krick-Pardo-Sombra}. 
We decided to use the cruder bound on $h(\delta)$ in Lemma~\ref{delta-lemma}
since it is independent of $n$.
\end{remark}

We now want to show that Koll\'ar's degree bound over fields
entails a similar bound for B\'ezout identities over rings of integers.

\begin{prop}\label{Bezout-Prop}
Suppose that $R={\cal O}_{F}$.
If $1\in (f_1,\dots,f_n)$, then there exist $g_1,\dots,g_n\in R[X]$ with
$$1=f_1h_1+\cdots+f_nh_n$$
and
$$\deg h_i \leq [F:\Q]\cdot (3d)^{O(N^2)} \bigl(h_1(f_1,\dots,f_n)+1\bigr)$$
for all $i=1,\dots,n$.
\end{prop}

Before we begin with the proof, we state an elementary lemma
whose proof is left to the reader:

\begin{lemma}
Let $U=(U_1,\dots,U_n)$, $V=(V_1,\dots,V_n)$ be tuples of pairwise
distinct indeterminates over $\Z$, and $e\geq 1$ an integer. There exist
polynomials $$g_1^{(e)}(U,V),\dots,g_n^{(e)}(U,V)$$
with non-negative integer coefficients such that
\begin{equation} \label{GenPower}
\bigl(1+U_1V_1+\cdots+U_nV_n\bigr)^e =
1+g_1^{(e)}(U,V)U_1+\cdots+g_n^{(e)}(U,V)U_n
\end{equation}
and $\deg_U g_j^{(e)}=e-1$, $\deg_V g_j^{(e)}=e$.
\end{lemma}

We first show a local analogue of Proposition~\ref{Bezout-Prop}:

\begin{lemma}\label{Local-Height}
Suppose that $R=({\cal O}_{F})_{\frak p}$, where $\frak p\neq 0$ is a prime 
ideal of ${\cal O}_{F}$. If $1\in (f_1,\dots,f_n)$, then
$$1=f_1h_1+\cdots+f_nh_n$$
for some $h_1,\dots,h_n\in R[X]$ of degree at most
\begin{equation}\label{hdeg}
[F:\Q]\cdot (3d)^{O(N^2)} 
\bigl(h(f_1,\dots,f_n)+1\bigr)/\log p.
\end{equation}
Here $p$ is the unique prime number such that $\Z\cap\frak p=p\Z$.
\end{lemma}
\begin{proof}
Suppose $1\in (f_1,\dots,f_n)$. Then $1\in (f_1,\dots,f_n)F[X]$,
hence by Lemma~\ref{delta-lemma} there exist
$g_1,\dots,g_n\in R[X]$ of degree at most $(3d)^N$ and a non-zero
$\delta\in R$ of height at most  
$\bigl(2(d+1)\bigr)^{O(N^2)}(h+1)$
such that
\begin{equation}\label{B1}
\delta = f_1g_1+\cdots+f_ng_n.
\end{equation}
Here and below $h=h(f_1,\dots,f_n)$. If $\delta$ is a unit in $R$, then
$$1=f_1(g_1/\delta)+\cdots+f_n(g_n/\delta)$$
is a B\'ezout identity for $f_1,\dots,f_n$ in $R[X]$, and $h_i:=g_i/\delta$,
$i=1,\dots,n$, have the required properties. Suppose that $\delta$ is not
a unit, so $e=v_{\frak p}(\delta)\geq 1$.
We have $1\in (\bar{f_1},\dots,\bar{f_n})$,
where $\bar{f}$ is the canonical image of $f\in R[X]$ in $(R/{\frak p}R)[X]$.
By Koll\'ar's theorem \cite{Kollar} applied to the field $R/{\frak p}R$, 
there exist $r_1,\dots,r_n\in R[X]$
with
$$1-(r_1f_1+\cdots+r_nf_n)\in{\frak p}R[X]$$
and $\deg r_j\leq (3d)^N$ for all $j=1,\dots,n$. Specializing the
$U_i$'s to $(f_1,\dots,f_n)$ and the $V_i$'s to $(-r_1,\dots,-r_n)$ in
\eqref{GenPower} gives 
$s_1,\dots,s_n\in R[X]$ and $s\in {\frak p}^e R[X]$ 
such that
\begin{equation}\label{B2}
1-(f_1s_1+\cdots+f_ns_n)=s.
\end{equation}
We have $\deg s_j \leq  e\bigl(d+(3d)^N\bigr)-d$ for all $j$, hence $\deg s\leq e\bigl(d+(3d)^N\bigr)$.
From \eqref{B1} and \eqref{B2} we get
$$1=f_1s_1+\cdots+f_ns_n+s=f_1h_1+\cdots+f_nh_n$$
with $h_j=s_j+(s/\delta)g_j\in R[X]$. We have
$$\deg (sg_j) \leq e\bigl(d+(3d)^N\bigr)+(3d)^N \leq 
e(3d)^{N+1},$$
and since $\bigl(2(d+1)\bigr)^{O(N^2+1)}=(3d)^{O(N^2)}$ 
we get
$$e\cdot \log p \leq [F:\Q]\cdot h(\delta)
= [F:\Q]\cdot (3d)^{O(N^2)} (h+1)$$
by the
remarks following Lemma~\ref{hv}, for $N>0$, $d>0$. 
It follows that $\deg h_j$ is bounded from
above by \eqref{hdeg}, for $j=1,\dots,n$.
\end{proof}

Now suppose that $R={\cal O}_F$, and assume that $1\in (f_1,\dots,f_n)$.
Hence by Lemma~\ref{delta-lemma} there are
$g_1,\dots,g_n\in R[X]$ of degree at most $(3d)^N$ and a non-zero
$\delta\in R$ of height at most  
$\bigl(2(d+1)\bigr)^{O(N^2+1)}(h+1)$
such that
$$\delta = f_1g_1+\cdots+f_ng_n.$$
For every prime ideal $\frak p$ of $R$ containing
$\delta$, and $p$ the prime number generating the ideal $\Z\cap\frak p$, 
we find $h_1^{(\frak p)},\dots,h_n^{(\frak p)}\in R[X]$ 
of degree bounded by \eqref{hdeg} as well as
$\delta^{(\frak p)}\in R\setminus\frak p$ such that
$$\delta^{(\frak p)}=f_1h_1^{(\frak p)}+\cdots+f_nh_n^{(\frak p)}.$$
Let ${\frak p}_1,\dots,{\frak p}_K$ be the pairwise distinct 
prime ideals of $R$ containing
$\delta$. Then there exist $a,a_1,\dots,a_K\in R$ such that
$$1=a\delta+a_1\delta^{({\frak p}_1)}+\cdots+a_K\delta^{({\frak p}_K)}.$$
Hence, letting $h_j=ag_j+a_1h_j^{({\frak p}_1)}+\cdots+a_Kh_j^{({\frak p}_K)}
\in R[X]$ for $j=1,\dots,n$, we get
\begin{multline*}
f_1h_1+\cdots+f_nh_n=a(f_1g_1+\cdots+f_ng_n)+\sum_{k=1}^K
a_k \bigl(f_1h_1^{({\frak p}_k)}+\cdots+f_nh_n^{({\frak p}_k)}\bigr)=\\
a\delta+\sum_{k=1}^K a_k\delta^{({\frak p}_k)} = 1.
\end{multline*}
From this Proposition~\ref{Bezout-Prop} follows. \qed

\subsection*{Ideal membership}
In the following we 
let $R={\cal O}_F$ for a number field $F$.
Let $A$ be an $m\times n$-matrix with entries in $R[X]$ and $b\in
\bigl(R[X]\bigr)^m$ a column vector.

\begin{theorem}\label{MainTh}
If the system
$Ay=b$
has a solution in $D=R[X]$, then it has such a solution of degree at most
\begin{equation}\label{degy}
[F:\Q]\cdot C_2\cdot\bigl(2m\deg(A,b)\bigr)^{2^{O(N^2)}}\cdot\bigl(h(A,b)+1\bigr).
\end{equation}
Here the constant $C_2$ depends only on $F$.
\end{theorem}
\begin{proof}
Put $d=\deg (A,b)$ and
$h=h(A,b)$.
By Corollary~\ref{HeightCor} there exist generators $z^{(1)},\dots,z^{(K)}$
for the $D$-module of solutions to the system of homogeneous linear equations
$(A,-b)z=0$, where $z$ is a vector of $n+1$ unknowns $z_1,\dots,z_{n+1}$,
with
\begin{align}
\deg\bigl(z^{(k)}\bigr) & = (2md)^{2^{O(N^2)}}, \label{d-bound}\\
h\bigl(z^{(k)}\bigr) & = C_2\cdot \bigl(2m(d+1)\bigr)^{O(N^2+1)}(h+1)\label{h1-bound}
\end{align}
for all $k=1,\dots,K$. The constant $C_2$ only depends on the number field $F$.
For each $k$ let $z_{n+1}^{(k)}\in R[X]$ be
the last component of $z^{(k)}$. Clearly, $Ay=b$ is solvable
in $R[X]$ if and only if $1\in\bigl(z_{n+1}^{(1)},\dots,z_{n+1}^{(K)}\bigr)$.
Moreover,
if $h_1,\dots,h_K$ are elements of $R[X]$ such that
$$1=h_1z_{n+1}^{(1)}+\cdots+h_Kz_{n+1}^{(K)},$$ then 
$y\in\bigl(R[X])^n$ with
$$\begin{bmatrix} y\\ 1\end{bmatrix}= h_1z^{(1)}+\cdots+h_Kz^{(K)}$$ 
is a solution to $Ay=b$. By Proposition~\ref{Bezout-Prop} we find such
$h_1,\dots,h_K$ with
$$\deg(h_k) \leq [F:\Q]\cdot \bigl(3\max_l\deg\bigl(z^{(l)}\bigr)\bigr)^{O(N^2)}
\bigl(\max_l h\bigl(z^{(l)}\bigr)+1\bigr)$$
for all $k$. 
From \eqref{d-bound}  we get
$$\bigl(3\max_l\deg\bigl(z^{(l)}\bigr)\bigr)^{O(N^2)}=
(2md)^{2^{O(N^2)}},$$
and using \eqref{h1-bound}
$$\bigl(\max_l h\bigl(z^{(l)}\bigr)+1\bigr) = 
C_2\cdot \bigl(2m(d+1)\bigr)^{2^{O(N^2+1)}}(h+1).$$
Hence the  vector $y$ has degree at
most \eqref{degy} as claimed.
\end{proof}

For $m=1$ the previous theorem yields:

\begin{cor}\label{MainThCor}
Let $f_0,f_1,\dots,f_n\in R[X]$, and put $d=\deg(f_1,\dots,f_n)$,
$h=h(f_1,\dots,f_n)$. If $f_0\in (f_1,\dots,f_n)$, then there
exist $g_1,\dots,g_n\in R[X]$ with $$f_0=g_1f_1+\cdots+g_nf_n$$ and
$$\deg(g_1,\dots,g_n) \leq [F:\Q]\cdot C_2\cdot (2d)^{2^{O(N^2)}}.$$
The constant $C_2$ depends only on $F$. \qed
\end{cor}

The doubly exponential degree bound on the solutions $y$ in 
Theorem~\ref{MainTh} implies a doubly exponential bound on $h(y)$. See
\cite{OLeary-Vaaler} for good bounds on the height of solutions to linear
equations over $R$.

Using the criterion for primeness of ideals in $\Z[X]$ given in
Corollary~\ref{PrimeCor}, we get:

\begin{cor}
One can test elementary recursively whether finitely many given polynomials
from $\Z[X]$ generate a prime ideal $I$.
\end{cor}
\begin{proof}
It is well-known that the conditions ``$I\Q[X]$ prime'' and ``$I\bmod p\subseteq\F_p[X]$ prime''
(for a prime $p$) can be tested elementary
recursively (\cite{Hermann}, \cite{Seidenberg1}). 
The condition ``$I\Q[X]\cap\Z[X]=I$''
may be tested elementary recursively using 
Proposition~\ref{ModuleOps} and Corollary~\ref{MainThCor}.
\end{proof}

See also \cite{GTZ} for an algorithm to test primeness of ideals in $\Z[X]$,
which is however not even obviously primitive recursive.

Similarly to Theorem~\ref{MainTh}, 
using Proposition~\ref{HeightProp} and Lemma~\ref{Local-Height}
in place of Corollary~\ref{HeightCor} and Proposition~\ref{Bezout-Prop}, 
respectively, one shows:

\begin{theorem}
Let $\frak p$ be a non-zero prime ideal of $R$ and $p$ the unique 
prime number with
$\Z\cap{\frak p}=p\Z$.
If the system
$Ay=b$
has a solution in $R_{\frak p}[X]$, then it has such a solution of degree at most
$$[F:\Q]\cdot C_2\cdot\bigl(2m\deg(A,b)\bigr)^{2^{O(N^2)}}\cdot\bigl(h(A,b)+1\bigr)/\log p.$$
Here $C_2$ is a constant which only depends on $F$.
\qed
\end{theorem}

The following example, promised in the introduction, shows that
the bounds established in the theorems above
necessarily have to depend not only on the degree, but also on the
coefficients of the polynomials involved.

\begin{example}
Let $p,d\in\Z$, $p>1$, $d\geq 1$. We have $1\in (1-pX,p^dX)\Z[X]$,
since
$$1=\bigl(1+pX+\cdots+p^{d-1}X^{d-1}\bigr)(1-pX)+X^{d-1}p^dX,$$
with the degrees of $1+pX+\cdots+p^{d-1}X^{d-1}$ and $X^{d-1}$ tending
to infinity, as $d\to\infty$. Considering everything mod $p^d$, we see
that $1-pX$ is a unit in $(\Z/p^d\Z)[X]$, indeed
$$1 \equiv \bigl(1+pX+\cdots+p^{d-1}X^{d-1}\bigr)(1-pX)\mod p^d,$$
so that if $1\equiv g(X)(1-pX)\bmod p^d$ for some $g(X)\in\Z[X]$, then
necessarily $g\equiv 1+pX+\cdots+p^{d-1}X^{d-1}\bmod p^d$. It follows
that if 
$$1=g(X)(1-pX)+h(X)p^dX\qquad\text{with $g,h\in\Z[X]$,}$$ 
then 
$\deg g,\deg h\geq d-1$.
Taking for $p$ a prime number and replacing $\Z$ by $\Z_{(p)}$, 
the same example works if we consider 
polynomials with coefficients in the ring $\Z_{(p)}$.
\end{example}

\subsection*{Final remarks}
From Theorem~\ref{MainTh} we obtain an effective 
reduction of the ideal membership problem
for $R[X]$, where $R={\cal O}_F$ for a number field $F$, to the
solvability of a (huge) system of linear equations over $R$.
As in the case of fields, this leads to a simple algorithm for deciding
ideal membership for $R[X]$. Certainly algorithms using 
Gr\"obner bases in $R[X]$ are much more effective in practice; it
remains to establish doubly exponential degree and height
bounds for the elements of 
Gr\"obner bases in $R[X]$. We plan to address this issue at a later point.

In \cite{Mayr}, Mayr shows that ideal membership problems
``$f_0\in (f_1,\dots,f_n)$''
with $f_0,f_1,\dots,f_n\in\Q[X]$ can be decided by an algorithm which uses space that grows
exponentially in the size of the input $f_0,\dots,f_n$.
Together with \cite{Mayr-Meyer} this establishes that
ideal membership in $\Q[X]$ is exponential-space complete.
The proof rests on an efficient parallel algorithm for computing the rank of
$m\times m$-matrices over $\Q$ in time $O(\log^2 m)$ and
the parallel computation thesis (``parallel time = sequential space'').
By \cite{Mayr-Meyer}, ideal membership in $\Z[X]$ is exponential-space hard.
Theorem~\ref{MainTh} (and the reduction given in
\cite{Mayr}) unfortunately only shows that ideal membership in $\Z[X]$ is
exponential-space complete {\em provided}\/ that solvability of systems of
linear equations over $\Z$ can be decided using logarithmic space.
However, this is even unknown for systems consisting of a single
equation of the form $1=ax+by$
($a,b\in\Z$), see \cite{GHR}.

\bibliographystyle{amsplain}
\bibliography{ideal2.bib}

\providecommand{\bysame}{\leavevmode\hbox to3em{\hrulefill}\thinspace}
\providecommand{\MR}{\relax\ifhmode\unskip\space\fi MR }
\providecommand{\MRhref}[2]{%
  \href{http://www.ams.org/mathscinet-getitem?mr=#1}{#2}
}
\providecommand{\href}[2]{#2}
\begin{thebibliography}{10}

\bibitem{coherent}
M.~Aschenbrenner, \emph{Bounds and definability in polynomial rings}, in
  preparation.

\bibitem{kron}
\bysame, \emph{Kronecker's problem}, in preparation.

\bibitem{Aschenbrenner-thesis}
\bysame, \emph{Ideal {M}embership in {P}olynomial {R}ings over the {I}ntegers},
  Ph.D. thesis, University of Illinois at Urbana-Champaign, 2001.

\bibitem{Ayoub}
C.~Ayoub, \emph{On constructing bases for ideals in polynomial rings over the
  integers}, J. Number Theory \textbf{17} (1983), no.~2, 204--225.

\bibitem{Baur}
W.~Baur, \emph{Rekursive {A}lgebren mit {K}ettenbedingungen}, Z. Math. Logik
  Grundlagen Math. \textbf{20} (1974), 37--46.

\bibitem{Becker-Weispfenning}
T.~Becker and V.~Weispfenning, \emph{Gr\"obner {B}ases}, Graduate Texts in
  Mathematics, vol. 141, Springer-Verlag, New York, 1993.

\bibitem{Berenstein-Yger}
C.~Berenstein and A.~Yger, \emph{Bounds for the degrees in the division
  problem}, Michigan Math. J. \textbf{37} (1990), 25--43.

\bibitem{BGR}
S.~Bosch, U.~G{\"u}ntzer, and R.~Remmert, \emph{Non-{A}rchimedean {A}nalysis.
  {A} {S}ystematic {A}pproach to {R}igid {A}nalytic {G}eometry}, Grundlehren
  der Mathematischen Wissenschaften, vol. 261, Springer-Verlag, Berlin, 1984.

\bibitem{Cohen}
H.~Cohen, \emph{A {C}ourse in {C}omputational {A}lgebraic {N}umber {T}heory},
  Graduate Texts in Mathematics, vol. 138, Springer-Verlag, Berlin, 1993.

\bibitem{DFGS}
A.~Dickenstein, N.~Fitchas, M.~Giusti, and C.~Sessa, \emph{The membership
  problem for unmixed polynomial ideals is solvable in single exponential
  time}, Discrete Appl. Math. \textbf{33} (1991), 73--94.

\bibitem{Edwards}
H.~M. Edwards, \emph{Kronecker's views on the foundations of mathematics}, The
  {H}istory of {M}odern {M}athematics, {V}ol.\ {I} \rom{(}{P}oughkeepsie, {NY},
  1989\rom{)}, Academic Press, Boston, MA, 1989, pp.~67--77.

\bibitem{Evans}
T.~Evans, \emph{Some connections between residual finiteness, finite
  embeddability and the word problem}, J. London Math. Soc. \rom{(2)}
  \textbf{1} (1969), 399--403.

\bibitem{Gallo-Mishra}
G.~Gallo and B.~Mishra, \emph{A solution to {K}ronecker's problem}, Appl.
  Algebra in Engrg. Comm. Comput. \textbf{5} (1994), no.~6, 343--370.

\bibitem{GTZ}
P.~Gianni, B.~Trager, and G.~Zacharias, \emph{Gr\"obner bases and primary
  decomposition of polynomial ideals}, J. Symbolic Comput. \textbf{6} (1988),
  no.~2-3, 149--167.

\bibitem{Gilmer}
R.~Gilmer, \emph{Multiplicative {I}deal {T}heory}, Queen's Papers in Pure and
  Applied Mathematics, vol.~12, Queen's University, Kingston, Ont., 1968.

\bibitem{Glaz}
S.~Glaz, \emph{Commutative {C}oherent {R}ings}, Lecture Notes in Math., vol.
  1371, Springer-Verlag, Berlin-Heidelberg-New York, 1989.

\bibitem{Greco-Salmon}
S.~Greco and P.~Salmon, \emph{Topics in $\m$-adic {T}opologies}, Ergebnisse der
  Mathe\-matik und ihrer Grenzgebiete, vol.~58, Springer-Verlag, New
  York-Berlin, 1971.

\bibitem{GHR}
R.~Greenlaw, H.~J. Hoover, and W.~L. Ruzzo, \emph{Limits to {P}arallel
  {C}omputation: {$P$}-{C}ompleteness {T}heory}, Oxford University Press,
  Oxford, 1995.

\bibitem{Hentzelt-Noether}
K.~Hentzelt and E.~Noether, \emph{Zur {T}heorie der {P}olynomideale und
  {R}esultanten}, Math. Ann. \textbf{88} (1923), 53--79.

\bibitem{Hermann}
G.~Hermann, \emph{{Die Frage der endlich vielen Schritte in der Theorie der
  Polynom\-ideale}}, Math. Ann. \textbf{95} (1926), 736--788.

\bibitem{KRK}
A.~Kandri-Rody and D.~Kapur, \emph{Computing a {G}r\"obner basis of a
  polynomial ideal over a {E}uclidean domain}, J. Symbolic Comput. \textbf{6}
  (1988), no.~1, 37--57.

\bibitem{Kollar}
J.~Koll\'ar, \emph{Sharp effective {N}ullstellensatz}, J. Amer. Math. Soc.
  \textbf{1} (1988), 963--975.

\bibitem{Krick-Pardo-2}
T.~Krick and L.~M. Pardo, \emph{A computational method for {D}iophantine
  approximation}, Algorithms in {A}lgebraic {G}eometry and {A}pplications
  \rom{(}{S}antander, 1994\rom{)}, Progr. Math., vol. 143, Birkh\"auser, Basel,
  1996, pp.~193--253.

\bibitem{Krick-Pardo-Sombra}
T.~Krick, L.~M. Pardo, and M.~Sombra, \emph{Sharp estimates for the arithmetic
  {N}ullstellensatz}, Duke Math. J. \textbf{109} (2001), no.~3, 521--598.

\bibitem{Lang-Diophantine}
S.~Lang, \emph{Fundamentals of {D}iophantine {G}eometry}, Springer-Verlag, New
  York, 1983.

\bibitem{Lang}
\bysame, \emph{Algebraic {N}umber {T}heory}, 2nd ed., Graduate Texts in
  Mathematics, vol. 110, Springer-Verlag, New York, 1994.

\bibitem{Mayr}
E.~Mayr, \emph{Membership in polynomial ideals over {${\cal Q}$} is exponential
  space complete}, STACS 89 (Paderborn, 1989), Lecture Notes in Comput. Sci.,
  vol. 349, Springer, Berlin, 1989, pp.~400--406.

\bibitem{Mayr-Meyer}
E.~Mayr and A.~Meyer, \emph{The complexity of the word problems for commutative
  semigroups and polynomial ideals}, Adv. Math. \textbf{46} (1982), no.~3,
  305--329.

\bibitem{Moreno-Socias}
G.~Moreno~Soc{\'{\i}}as, \emph{Length of polynomial ascending chains and
  primitive recursiveness}, Math. Scand. \textbf{71} (1992), no.~2, 181--205.

\bibitem{OLeary-Vaaler}
R.~O'Leary and J.~Vaaler, \emph{Small solutions to inhomogeneous linear
  equations over number fields}, Trans. Amer. Math. Soc. \textbf{336} (1993),
  no.~2, 915--931.

\bibitem{Renschuch}
B.~Renschuch, \emph{Beitr\"age zur konstruktiven {T}heorie der {P}olynomideale.
  {XVII}/1. {Z}ur {H}ent\-zelt/{N}oether/{H}ermannschen {T}heorie der endlich
  vielen {S}chritte}, Wiss. Z. P\"adagog. Hochsch. ``Karl Liebknecht'' Potsdam
  \textbf{24} (1980), no.~1, 87--99.

\bibitem{Richman}
F.~Richman, \emph{Constructive aspects of {N}oetherian rings}, Proc. Amer.
  Math. Soc. \textbf{44} (1974), 436--441.

\bibitem{Roy-Thunder}
D.~Roy and J.~L. Thunder, \emph{Bases of number fields with small height},
  Rocky Mountain J. Math. \textbf{26} (1996), no.~3, 1089--1098.

\bibitem{Schmidt-Goettsch}
K.~Schmidt-G\"ottsch, \emph{Polynomial bounds in polynomial rings over fields},
  J. Algebra \textbf{125} (1989), 164--180.

\bibitem{Seidenberg1}
A.~Seidenberg, \emph{Constructions in algebra}, Trans. Amer. Math. Soc.
  \textbf{197} (1974), 273--313.

\bibitem{Seidenberg3}
\bysame, \emph{What is {N}oetherian?}, Rend. Sem. Mat. Fis. Milano \textbf{44}
  (1974), 55--61.

\bibitem{Simmons}
H.~Simmons, \emph{The solution of a decision problem for several classes of
  rings}, Pacific J. Math. \textbf{34} (1970), 547--557.

\bibitem{Sombra}
M.~Sombra, \emph{A sparse effective {N}ullstellensatz}, Adv. in Appl. Math.
  \textbf{22} (1999), 271--295.

\bibitem{Wainer}
S.~S. Wainer, \emph{A classification of the ordinal recursive functions}, Arch.
  Math. Logik Grundlagenforsch. \textbf{13} (1970), 136--153.

\end{thebibliography}

\end{document}